\documentclass[runningheads]{llncs}
\usepackage{mathtools} %
\usepackage{amssymb} %

\usepackage{amsthm} %
\usepackage{stmaryrd} %
\usepackage{relsize} %
\usepackage{microtype} %
\usepackage{multicol} %
\usepackage{csquotes} %
\usepackage{xspace}
\usepackage{enumitem}
\usepackage[nocompress]{cite} %
  \usepackage{graphicx} %
  \usepackage{subcaption} %
  \usepackage[dvipsnames]{xcolor} %

    \usepackage{tikz} %
    \usetikzlibrary{
      cd, %
      petri, %
      backgrounds, %
      arrows, %
      positioning, %
      decorations.markings, %
      calc,  %
      fit, %
    }

\usepackage{tabularx, cellspace} %
\usepackage{hyperref}
\usepackage[capitalize]{cleveref}
\usepackage{subfiles}

  \newcounter{theoremUnified} %
  \numberwithin{theoremUnified}{section} %
  \numberwithin{theoremUnified}{section} %

  \newtheoremstyle{plainStyle} %
  {2mm} %
  {2mm} %
  {} %
  {} %
  {\bfseries} %
  {.} %
  {.5em} %
  {} %

  \newtheoremstyle{italicStyle} %
  {2mm} %
  {2mm} %
  {\itshape} %
  {} %
  {\bfseries} %
  {.} %
  {.5em} %
  {} %

\tikzstyle{place}=
[circle,thick,draw=blue!75,fill=blue!20,minimum size=6mm]
\tikzstyle{transition}=
[rectangle,thick,draw=black!75,fill=black!20,minimum size=4mm]

\newcommand{\cp}{\fatsemi} %

\newcommand{\Obj}[1]{\operatorname{Obj} \, #1} %
\newcommand{\GObj}[1]{\operatorname{GenObj} \, #1} %
\newcommand{\GMor}[1]{\operatorname{GenMor} \, #1} %

\newcommand{\Hom}[3]{\operatorname{Hom}_{\,#1}\left[#2,#3\right]} %
\newcommand{\Id}[1]{id_{#1}} %

\newcommand{\CategoryC}{\mathcal{C}}
\newcommand{\CategoryD}{\mathcal{D}}

\newcommand{\TermCat}{\textbf{1}} %
\newcommand{\PetriS}[1]{\textbf{Petri}^{#1}} %
  \newcommand{\PetriSetS}{\PetriS{\SetS}} %
  \newcommand{\PetriSpan}{\PetriS{\Span}} %
  \newcommand{\PetriTerm}{\PetriS{\TermCat}} %

\newcommand{\Fun}[1]{{#1}^\sharp} %
\newcommand{\NetSem}[1]{\left( #1, \Fun{#1}\right)} %
  \newcommand{\NetSemTerm}[1]{\NetSem{#1}} %

\newcommand{\Semantics}{\mathcal{S}} %

\newcommand{\Free}[1]{\mathfrak{F}\left(#1\right)} %
\newcommand{\UnFree}[1]{\mathfrak{U}\left(#1\right)} %

\newcommand{\Grothendieck}[1]{\textstyle\int{#1}} %
\newcommand{\GrothendieckS}[1]{\Grothendieck{\Fun{#1}}} %

\newcommand{\EmbSetS}{\texttt{emb}_\SetS} %
\newcommand{\EmbSpan}{\texttt{emb}_\Span} %
\newcommand{\ForSetS}{\texttt{for}_\SetS} %
\newcommand{\ForSpan}{\texttt{for}_\Span} %

\newcommand{\harp}[1]{\mathpalette\harpoonvec{#1}}
\newcommand{\harpvecsign}{\scriptscriptstyle\rightharpoonup}
\newcommand{\harpoonvec}[2]{%
  \ifx\displaystyle#1\doalign{$\harpvecsign$}{#1#2}\fi
  \ifx\textstyle#1\doalign{$\harpvecsign$}{#1#2}\fi
  \ifx\scriptstyle#1\doalign{\scalebox{.6}[.9]{$\harpvecsign$}}{#1#2}\fi
  \ifx\scriptscriptstyle#1\doalign{\scalebox{.5}[.8]{$\harpvecsign$}}{#1#2}\fi
}
\newcommand{\doalign}[2]{%
 {\vbox{\offinterlineskip\ialign{\hfil##\hfil\cr#1\cr$#2$\cr}}}%
}

\newcommand{\LiftSetS}[1]{\harp{#1}} %
\newcommand{\LiftSpan}[1]{\widehat{#1}} %

\newcommand{\Set}{\textbf{Set}} %
\newcommand{\SetS}{{\Set_{*}}} %
\newcommand{\Rel}{\textbf{Rel}} %
\newcommand{\Span}{\textbf{Span}} %
\newcommand{\Cat}{\textbf{Cat}} %

\newcommand{\Tensor}{\otimes} %
\newcommand{\TensorUnit}{I} %

\newcommand{\Suchthat}[2]{\left\{#1 \: \middle\vert \: #2\right\}} %

\tikzset{ %
  oriented WD/.style={%
    every to/.style={
      out=0,in=180,draw
    },
    label/.style={
      font=\everymath\expandafter{\the\everymath\scriptstyle},
      inner sep=0pt,
      node distance=2pt and -2pt
    },
    semithick,
    node distance=1 and 1,
    decoration={
      markings, mark=at position \stringdecpos with \stringdec
    },
    ar/.style={
      postaction={decorate}
    },
    execute at begin picture={
      \tikzset{
        x=\bbx, y=\bby,
        every fit/.style={
          inner xsep=\bbx, inner ysep=\bby
        }
      }
    }
  },
  string decoration/.store in=\stringdec,
  string decoration={
    \arrow{stealth};
  },
  string decoration pos/.store in=\stringdecpos,
  string decoration pos=.7,
  bbx/.store in=\bbx,
  bbx = 1.5cm,
  bby/.store in=\bby,
  bby = 1.5ex,
  bb port sep/.store in=\bbportsep,
  bb port sep=1.5,
  bb port length/.store in=\bbportlen,
  bb port length=4pt,
  bb penetrate/.store in=\bbpenetrate,
  bb penetrate=0,
  bb min width/.store in=\bbminwidth,
  bb min width=1cm,
  bb rounded corners/.store in=\bbcorners,
  bb rounded corners=2pt,
  bb small/.style={
    bb port sep=1, 
    bb port length=2.5pt, 
    bbx=.4cm, bb min width=.4cm, 
    bby=.7ex
  },
  bb medium/.style={
    bb port sep=1, 
    bb port length=2.5pt, 
    bbx=.4cm, 
    bb min width=.4cm, 
    bby=.9ex
  },
  bb/.code 2 args={%
    \pgfmathsetlengthmacro{\bbheight}{\bbportsep * (max(#1,#2)+1) * \bby}
    \pgfkeysalso{
      draw,
      minimum height=\bbheight,
      minimum width=\bbminwidth,
      outer sep=0pt,
      rounded corners=\bbcorners,
      thick,
      prefix after command={
        \pgfextra{\let\fixname\tikzlastnode}
      },
      append after command={
        \pgfextra{
          \draw
          \ifnum #1=0
            {} 
          \else 
            foreach \i in {1,...,#1} {
              ($(\fixname.north west)!{\i/(#1+1)}!(\fixname.south west)$) +(-
            \bbportlen,0) 
            coordinate (\fixname_in\i) -- +(\bbpenetrate,0) coordinate (\fixname_in\i')
            }
          \fi 
          \ifnum 
            #2=0{} 
          \else 
            foreach \i in {1,...,#2} {
            ($(\fixname.north east)!{\i/(#2+1)}!(\fixname.south east)$) +(-
            \bbpenetrate,0) 
            coordinate (\fixname_out\i') -- +(\bbportlen,0) coordinate (\fixname_out\i)
            }
          \fi;
        }
      }
    }
  },
  bb name/.style={
    append after command={
      \pgfextra{
        \node[anchor=north] at (\fixname.north) {#1}
      ;}
    }
  }
}

\newcommand{\from}{\leftarrow}
\newcommand{\From}[1]{\xleftarrow{#1}}
\newcommand{\To}[1]{\xrightarrow{#1}}

\begin{document}
\title{A Categorical Semantics for Guarded Petri Nets}

\author{Fabrizio Genovese\inst{1}\orcidID{0000-0001-7792-1375} \and
David I. Spivak\inst{2}\orcidID{0000-0002-9326-5328}}

\authorrunning{F. Genovese et al.}

\institute{Statebox\\ \email{research@statebox.io}
\and
MIT\\ \email{dspivak@mit.edu}}
\maketitle   
\begin{abstract}
  We build on the correspondence between Petri nets and free symmetric 
  strict monoidal categories already investigated in the literature, and present 
  a categorical semantics for Petri nets with guards. This comes in two flavors: 
  Deterministic and with side-effects. Using the Grothendieck construction, we show 
  how the guard semantics can be internalized in the net itself.
\end{abstract}
\section{Introduction}
  \label{sec: introduction}
Category theory has been used to study Petri nets 
at least since the beginning of the nineties~\cite{Meseguer1990}. Throughout this time, 
the main effort in this direction of research consisted in showing 
how Petri nets can be thought of as presenting various flavors of
free monoidal categories~\cite{Meseguer1990, Sassone1995,Genovese2019b,Master2019} 
This idea has been very influential, 
successfully modeling the individual-token philosophy via 
process semantics.

On the other hand, shortly after Petri's first publications 
about the nets that carry his name~\cite{Petri2008} researchers started
investigating what happens when nets are enriched with 
new features. One of the most successful extensions of 
Petri nets is \emph{guarded (or coloured) nets}~\cite{Jensen2009}.
Modulo different flavors of modeling what boils 
down to be the same concept, a \emph{guarded net} is a 
Petri net with the following extra properties:
\begin{itemize}
  \item To each token is attached some ``attribute''. The kind 
  of attributes we can attach to tokens depends on the place 
  the token is in;
  \item Each arc is decorated with an expression, which modifies tokens' 
  attributes as they flow through the net;
  \item Each transition is decorated with a predicate and only fires on tokens whose attributes satisfy 
  the predicate.
\end{itemize}
At a fist glance, guarded nets allow for a more expressive 
form of modeling with respect to their unguarded counterparts, 
but as we will see shortly, this is not necessarily the case. Indeed, 
depending on the underlying theory from which properties, expressions, 
and predicates are drawn the gain in expressive power with respect 
to undecorated nets may be nil:  With a wise choice of 
underlying theory, coloured nets amount to be nothing more than 
syntactic sugar for standard nets, though of course the availability of such syntactic 
sugar can greatly simplify the modeling of complex processes 
using the Petri net formalism.

Recently there has been renewed interest in employing
Petri nets as the basis for a programming language~\cite{StateboxTeam2017}.
In this setting, the categorical correspondence between nets 
and symmetric monoidal categories has been of the utmost importance, 
single-token philosophy being considered necessary to make 
the programming language usable~\cite{StateboxTeam2019}. Clearly, extending nets 
with new features such as guards or timings is desirable to 
make the language more expressive.

In this work we try to unify these two longstanding directions 
of research -- the categorical approach to Petri nets 
and the study of guarded nets -- by showing how guarded nets can be modeled as ordinary
Petri nets with a particular flavor of semantics in the style of~\cite{Genovese2019}.

Importantly, we are able to define both a \emph{deterministic 
semantics} and a \emph{non-deterministic semantics} in 
our formalism. The first models the traditional notion of guards 
deterministically modifying data attached to tokens, while 
the second describes a setting where token data is 
modified depending on side effects.

Using the Grothendieck construction, we show 
how the guard semantics can be \emph{internalized} in the net 
itself, providing a categorical proof that in our model, guarded nets do not 
increase expressivity, as compared to traditional nets. 
This is a desired feature, since it means that many nice properties 
of nets such as termination or decidability of the reachability 
relation are preserved. It also shows that the core mathematical abstraction
in computer implementations of Petri nets need not be modified when offering users
the flexibility of guarded nets.

We save all proofs for the appendix, which starts on page~\pageref{appendix}.
\section{Guarded nets}
  \label{sec: guarded nets}
Having given an intuitive version of what a 
guarded net is, we now start modeling the concept 
formally. We will use the formalism developed in~\cite{Genovese2019},
of which we recall some core concepts.

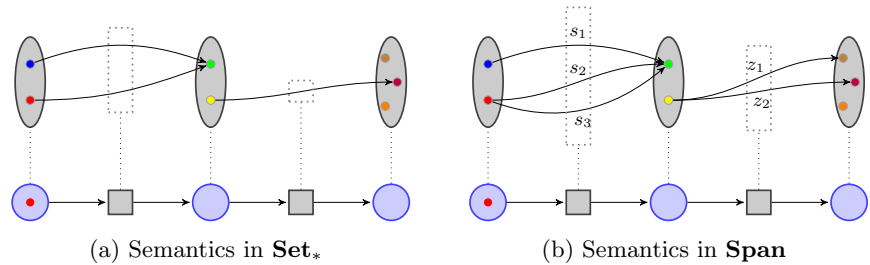
\begin{figure}[!ht]\centering
  \begin{subfigure}[t]{0.49\textwidth}\centering
    \scalebox{0.8}{
    \begin{tikzpicture}[node distance=1.3cm,>=stealth',bend angle=45,auto]
      \node [place,colored tokens={red},] (1a) at (0,0){};
      \node [transition] (2a) at (1.5,0)     {}
          edge [pre] (1a);
      \node [place,tokens=0] (3a) at (3,0)  {}
          edge [pre] (2a);
      \node [transition] (4a) at (4.5,0) {}
          edge [pre] (3a);
      \node [place,tokens=0] (5a) at (6,0) {}
          edge [pre] (4a);
      \draw[thick, draw=black!75, fill=black!20] (0,2) ellipse (0.25cm and 0.75cm);
      \draw[thick, draw=black!75, fill=black!20] (3,2) ellipse (0.25cm and 0.75cm);
      \draw[thick, draw=black!75, fill=black!20] (6,2) ellipse (0.25cm and 0.75cm);

      \node[circle,fill=red, draw=gray!75, inner sep=0pt,minimum size=4pt] (red) at (0,1.7) {};
      \node[circle,fill=blue, draw=gray!75, inner sep=0pt,minimum size=4pt] (blue) at (0,2.3) {};

      \node[circle,fill=yellow, draw=gray!75, inner sep=0pt,minimum size=4pt] (yellow) at (3,1.7) {};
      \node[circle,fill=green, draw=gray!75, inner sep=0pt,minimum size=4pt] (green) at (3,2.3) {};

      \node[circle,fill=brown, draw=gray!75, inner sep=0pt,minimum size=4pt] (brown) at (5.9,2.4) {};
      \node[circle,fill=purple, draw=gray!75, inner sep=0pt,minimum size=4pt] (purple) at (6.1,2) {};
      \node[circle,fill=orange, draw=gray!75, inner sep=0pt,minimum size=4pt] (orange) at (5.9,1.6) {};

      \draw[->, out=20, in=160] (blue) to (green);
      \draw[->, out=0, in=200] (red) to (green);
      \draw[->, out=0, in=180] (yellow) to (purple);
      
      \draw[dotted, -]  (1a.north) -- (0,1.25);
      \draw[dotted, -]  (3a.north) -- (3,1.25);
      \draw[dotted, -]  (5a.north) -- (6,1.25);

      \node[transition,dotted,fill=none, draw=black!50, inner sep=0pt,minimum height=40pt] (left) at (1.5,2.2) {};
      \node[transition,dotted,fill=none, draw=black!50, inner sep=0pt,minimum height=10pt] (right) at (4.5,1.85) {};

      \draw[dotted, -]  (2a.north) -- (left.south);
      \draw[dotted, -]  (4a.north) -- (right.south);

    \end{tikzpicture}}
    \caption{Semantics in $\SetS$}\label{fig: semantics in SetS}
  \end{subfigure}
  \begin{subfigure}[t]{0.49\textwidth}\centering
    \scalebox{0.8}{
    \begin{tikzpicture}[node distance=1.3cm,>=stealth',bend angle=45,auto]
      \node [place,colored tokens={red},] (1a) at (0,0){};
      \node [transition] (2a) at (1.5,0)     {}
          edge [pre] (1a);
      \node [place,tokens=0] (3a) at (3,0)  {}
          edge [pre] (2a);
      \node [transition] (4a) at (4.5,0) {}
          edge [pre] (3a);
      \node [place,tokens=0] (5a) at (6,0) {}
          edge [pre] (4a);
      \draw[thick, draw=black!75, fill=black!20] (0,2) ellipse (0.25cm and 0.75cm);
      \draw[thick, draw=black!75, fill=black!20] (3,2) ellipse (0.25cm and 0.75cm);
      \draw[thick, draw=black!75, fill=black!20] (6,2) ellipse (0.25cm and 0.75cm);

      \node[circle,fill=red, draw=gray!75, inner sep=0pt,minimum size=4pt] (red) at (0,1.7) {};
      \node[circle,fill=blue, draw=gray!75, inner sep=0pt,minimum size=4pt] (blue) at (0,2.3) {};

      \node[circle,fill=yellow, draw=gray!75, inner sep=0pt,minimum size=4pt] (yellow) at (3,1.7) {};
      \node[circle,fill=green, draw=gray!75, inner sep=0pt,minimum size=4pt] (green) at (3,2.3) {};

      \node[circle,fill=brown, draw=gray!75, inner sep=0pt,minimum size=4pt] (brown) at (5.9,2.4) {};
      \node[circle,fill=purple, draw=gray!75, inner sep=0pt,minimum size=4pt] (purple) at (6.1,2) {};
      \node[circle,fill=orange, draw=gray!75, inner sep=0pt,minimum size=4pt] (orange) at (5.9,1.6) {};

      \draw[->, out=20, in=160] (blue) to node[midway, above] {$s_1$} (green);
      \draw[->, out=0, in=180] (red) to node[midway, above] {$s_2$} (green);
      \draw[->, out=-20, in=220] (red) to node[midway, below] {$s_3$} (green);
      \draw[->, out=0, in=180] (yellow) to node[midway, below] {$z_2$} (purple);
      \draw[->, out=0, in=180] (yellow) to node[midway, above] {$z_1$} (brown);
      
      \draw[dotted, -]  (1a.north) -- (0,1.25);
      \draw[dotted, -]  (3a.north) -- (3,1.25);
      \draw[dotted, -]  (5a.north) -- (6,1.25);

      \node[transition,dotted,fill=none, draw=black!50, inner sep=0pt,minimum height=65pt] (left) at (1.5,2.1) {};
      \node[transition,dotted,fill=none, draw=black!50, inner sep=0pt,minimum height=40pt] (right) at (4.5,1.9) {};

      \draw[dotted, -]  (2a.north) -- (left.south);
      \draw[dotted, -]  (4a.north) -- (right.south);
    \end{tikzpicture}}
    \caption{Semantics in $\Span$}\label{fig: semantics in span}
  \end{subfigure}
  \caption{The same net (below), equipped with a partial function and span semantics, respectively (above).}
\end{figure}

We denote by $\Free{N}$ the free symmetric strict 
monoidal category associated to a Petri net $N$,
and with $\UnFree{\CategoryC}$ the petri net associated 
to the free symmetric strict monoidal category $\CategoryC$. We denote composition in diagrammatic order; i.e.\ given $f\colon c\to d$ and $g\colon d\to e$, we denote their composite by $(f\cp g)\colon c\to e$.
\begin{definition}\label{def: PetriS}
  Given a strict monoidal category $\Semantics$, a 
  \emph{Petri net with $\Semantics$-semantics}
   is a pair  $\NetSem{N}$, consisting of a Petri net $N$ and
   a strict monoidal functor
  \[\Fun{N}: \Free{N} \to \Semantics.\]
  A morphism $F: \NetSem{M} \to \NetSem{N}$
  is just a strict monoidal functor $F: \Free{M} \to \Free{N}$ such that 
  $\Fun{M} = F \cp \Fun{N}$.

  Nets equipped with $\Semantics$-semantics and their morphisms form a monoidal category
  denoted $\PetriS{\Semantics}$, with the monoidal structure arising from the
  product in $\Cat$.
\end{definition}
\begin{definition}
  We denote by $\SetS$ the category of sets and partial functions, and by $\Span$ 
  the 1-category of sets and spans, where isomorphic spans are identified. Both 
  these categories are symmetric monoidal. From now on, we will 
  work with the \emph{strictified} version of $\SetS$ and $\Span$, respectively.
\end{definition}
\begin{example}\label{ex.petri_one}
Let $\TermCat$ denote the terminal symmetric monoidal category. A Petri net with 
$\TermCat$-semantics is just a Petri net. Petri nets are in bijective correspondence with
free symmetric strict monoidal categories, so $\PetriTerm$ denotes the usual category 
of free symmetric strict monoidal categories and strict monoidal functors between them.
\end{example}
\begin{notation}
  Recall that a morphism $A\to B$ in $\Span$ consists of a set $S$ and a pair 
  of functions $A\from S\to B$. When we need to notationally extract this 
  data from $f$, we write
  \begin{equation*}
    A\From{f_1}S_f\To{f_2}B
  \end{equation*}
  We sometimes consider the span as a function $f\colon S_f\to A\times B$, 
  thus we may write $f(s)=(a,b)$ for $s\in S_f$ with $f_1(s)=a$ and $f_2(s)=b$.
\end{notation}
Perhaps unsurprisingly, $\SetS$ and $\Span$ will be the target 
semantics corresponding to two different flavors for 
our guards, with $\Span$ allowing for some form 
of nondeterminism -- expressed as the action of 
side-effects -- whereas $\SetS$ models a purely 
deterministic semantics. Expressing things formally:
\begin{definition}
  \label{def: guarded net}
  A \emph{guarded net} is an object of $\PetriSetS$.
  A \emph{guarded net with side effects} is an object 
  of $\PetriSpan$. A morphism of guarded nets (with side effects) 
  is a morphism in $\PetriSetS$ (resp.\ in $\PetriSpan$).
\end{definition}
\begin{remark}
  Although it doesn't affect our formalism by any means, 
  in practice the choice of semantics, both for $\SetS$ and $\Span$, 
  is limited by computational requirements: the places in a net are usually 
  sent to \emph{finite} sets, while transitions are usually sent 
  to computable functions and spans\footnote{
    A \emph{computable span} is one for which both legs are 
    computable functions.
  }, respectively. Such restrictions are necessary to make sure the 
  net is executable and to keep model checking decidable.
\end{remark}
Let us unroll the cryptic \cref{def: guarded net}, starting from the case $\PetriSetS$.
An object in $\PetriSetS$ is a net $N$ together with a strict monoidal functor $\Fun{N}\colon\Free{N} \to \SetS$.
It assigns to each place $p$ of $N$ -- corresponding to a generating object of 
$\Free{N}$ -- a set $\Fun{N}(p)$, representing all the possible colours a token in $p$
can assume. A transition $f\colon p\to p'$ -- corresponding to a generating morphism
of $\Free{N}$ -- gets sent to a partial function $\Fun{N}(f)\colon \Fun{N}(p)\to \Fun{N}(p')$, representing 
how token colours are transformed during firing.
Importantly, the fact that the functions in the semantics are 
\emph{partial} means that a transition may not be defined for 
tokens of certain colors. An example of this is the net 
in \cref{fig: semantics in SetS}, which is 
shown together with its semantics. Although reachability 
in the base net seems quite straightforward, we see that a token 
in the leftmost place will never reach the rightmost place, since the
rightmost transition is not defined on the tokens output by the leftmost 
one.

In the case of $\PetriSpan$ the intuition is similar. Objects are sent 
to sets, exactly as in $\PetriSetS$, but transitions are mapped to spans.
Spans can be understood as \emph{relations with witnesses}, provided by 
elements in the apex of the span. Practically, this means that each path from the span 
domain to its codomain is indexed by some element of the span apex, as it  
is shown in \cref{fig: semantics in span}. The presence of witnesses 
allows to consider different paths between the same elements. Moreover, 
an element in the domain can be sent to different elements in the codomain
via different paths. We interpret this as \emph{non-determinism}: The 
firing of the transition is not only a matter of the tokens input and output, it 
also includes the path chosen, which we interpret as having 
side-effects that are interpreted outside of our model.
As one can see, in both \cref{fig: semantics in SetS,fig: semantics in span} the composition 
of paths is the empty function (resp.\ span). Seeing things from 
a reachability point of view, the process given by firing the left transition 
and then the right will never occur. Let us make this intuition precise:
\begin{definition}
  Given a guarded Petri net (with side effects) $\NetSem{N}$, a \emph{marking}
  for $\NetSem{N}$ is a pair $(X, x)$ where $X$ is an object of $\Free{N}$ and
  $x \in \Fun{N}X$. We say that a marking $(Y,y)$ is \emph{reachable} from $(X,x)$
  if there is a morphism $f: X \to Y$ in $\Free{N}$ such that $\Fun{N}f(x) = y$.
\end{definition}
The goal we will pursue 
in the next section will be to internalize the guard semantics
in the free category $\Free{N}$ associated to a net.
\section{Internalizing guards}
  \label{sec: internalizing guards}
By ``internalizing the semantics of a guarded net $N$ in $\Free{N}$''
we mean \emph{obtaining an unguarded net $M$ such that $\Free{M}$
represents all the possible runs of $N$}. For readers familiars with
coloured Petri nets, this corresponds to the claim that reachability in
a coloured net is equivalent to reachability in a suitably constructed
``standard'' net~\cite{Jensen2009}.

Since our point of view is process-theoretic, and we are working
with symmetric strict monoidal categories and functors, such internalization
must be built categorically. The main tool we will use is the
\emph{Grothendieck construction}~\cite{MacLane1992}, which in our context we will 
specialize to functors to $\SetS$ and $\Span$, respectively.
\begin{definition}\label{def.grothendieck_setS_case}
  Let $\NetSem{M}\in\PetriSetS$ be a guarded net. We define 
  its \emph{internalization}, denoted $\Grothendieck{\Fun{M}}$, as the following category:
  \begin{itemize}
    \item The objects of $\GrothendieckS{M}$ are pairs 
    $(X, x)$, where $X$ is an object of $\Free{M}$ and $x$ is an 
    element of $\Fun{M}X$. Concisely:
    \begin{equation*}
      \Obj{\GrothendieckS{M}} := 
        \Suchthat{(X,x)}{(X \in \Obj{\Free{M}}) \wedge (x \in \Fun{M}X)}.
    \end{equation*}
    \item A morphism from $(X,x)$ to 
    $(Y,y)$ in  $\Grothendieck{\Fun{M}}$ is a morphism $f\colon X \to Y$ in $\Free{M}$ such that 
    $x$ is sent to $y$ via $\Fun{M}f$. Concisely:
    \begin{equation*}
      \Hom{\GrothendieckS{M}}{(X,x)}{(Y,y)} :=
        \Suchthat{f}{(f \in \Hom{\Free{M}}{X}{Y}) \wedge (\Fun{M}f(x) = y)}.
    \end{equation*}
  \end{itemize}
\end{definition}
It is worth giving some intuition of what the Grothendieck construction
does in our context. It basically makes a place for each element of the set we send a place 
to, and makes a transition for each path between these elements, 
as shown below:
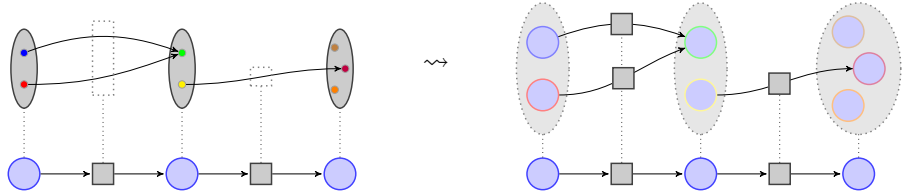
\begin{equation*}
  \scalebox{0.7}{
    \begin{tikzpicture}[node distance=1.3cm,>=stealth',bend angle=45,baseline=2cm]
      \node [place, tokens=0,] (1a) at (0,0){};
      \node [transition] (2a) at (1.5,0)     {}
          edge [pre] (1a);
      \node [place,tokens=0] (3a) at (3,0)  {}
          edge [pre] (2a);
      \node [transition] (4a) at (4.5,0) {}
          edge [pre] (3a);
      \node [place,tokens=0] (5a) at (6,0) {}
          edge [pre] (4a);
      \draw[thick, draw=black!75, fill=black!20] (0,2) ellipse (0.25cm and 0.75cm);
      \draw[thick, draw=black!75, fill=black!20] (3,2) ellipse (0.25cm and 0.75cm);
      \draw[thick, draw=black!75, fill=black!20] (6,2) ellipse (0.25cm and 0.75cm);

      \node[circle,fill=red, draw=gray!75, inner sep=0pt,minimum size=4pt] (red) at (0,1.7) {};
      \node[circle,fill=blue, draw=gray!75, inner sep=0pt,minimum size=4pt] (blue) at (0,2.3) {};

      \node[circle,fill=yellow, draw=gray!75, inner sep=0pt,minimum size=4pt] (yellow) at (3,1.7) {};
      \node[circle,fill=green, draw=gray!75, inner sep=0pt,minimum size=4pt] (green) at (3,2.3) {};

      \node[circle,fill=brown, draw=gray!75, inner sep=0pt,minimum size=4pt] (brown) at (5.9,2.4) {};
      \node[circle,fill=purple, draw=gray!75, inner sep=0pt,minimum size=4pt] (purple) at (6.1,2) {};
      \node[circle,fill=orange, draw=gray!75, inner sep=0pt,minimum size=4pt] (orange) at (5.9,1.6) {};

      \draw[->, out=20, in=160] (blue) to (green);
      \draw[->, out=0, in=200] (red) to (green);
      \draw[->, out=0, in=180] (yellow) to (purple);
      
      \draw[dotted, -]  (1a.north) -- (0,1.25);
      \draw[dotted, -]  (3a.north) -- (3,1.25);
      \draw[dotted, -]  (5a.north) -- (6,1.25);

      \node[transition,dotted,fill=none, draw=black!50, inner sep=0pt,minimum height=40pt] (left) at (1.5,2.2) {};
      \node[transition,dotted,fill=none, draw=black!50, inner sep=0pt,minimum height=10pt] (right) at (4.5,1.85) {};

      \draw[dotted, -]  (2a.north) -- (left.south);
      \draw[dotted, -]  (4a.north) -- (right.south);

    \end{tikzpicture}
  }
  \qquad \leadsto \qquad
  \scalebox{0.7}{
    \begin{tikzpicture}[node distance=1.3cm,>=stealth',bend angle=45, baseline=2cm]
      \node [place, tokens=0,] (1a) at (0,0){};
      \node [transition] (2a) at (1.5,0)     {}
          edge [pre] (1a);
      \node [place,tokens=0] (3a) at (3,0)  {}
          edge [pre] (2a);
      \node [transition] (4a) at (4.5,0) {}
          edge [pre] (3a);
      \node [place,tokens=0] (5a) at (6,0) {}
          edge [pre] (4a);
      \draw[thick, dotted, draw=black!50, fill=black!10] (0,2) ellipse (0.5cm and 1.25cm);
      \draw[thick, dotted, draw=black!50, fill=black!10] (3,2) ellipse (0.5cm and 1.25cm);
      \draw[thick, dotted, draw=black!50, fill=black!10] (6,2) ellipse (0.8cm and 1.25cm);

      \node[place, draw=red!50, tokens=0] (red) at (0,1.5) {};
      \node[place, draw=blue!50, tokens=0] (blue) at (0,2.5) {};

      \node[place, draw=yellow!50, tokens=0] (yellow) at (3,1.5) {};
      \node[place, draw=green!50, tokens=0] (green) at (3,2.5) {};

      \node[place, draw=brown!50, tokens=0] (brown) at (5.8,2.7) {};
      \node[place, draw=purple!50, tokens=0] (purple) at (6.2,2) {};
      \node[place, draw=orange!50, tokens=0] (orange) at (5.8,1.3) {};

      \draw[dotted, -]  (2a.north) -- (1.5,3);
      \draw[dotted, -]  (4a.north) -- (4.5,1.8);

      \draw[->, out=20, in=160] (blue) to node[transition] {} (green);
      \draw[->, out=0, in=200] (red) to node[transition] {} (green);
      \draw[->, out=0, in=180] (yellow) to node[pos=0.455, transition] {} (purple);
      
      \draw[dotted, -]  (1a.north) -- (0,0.75);
      \draw[dotted, -]  (3a.north) -- (3,0.75);
      \draw[dotted, -]  (5a.north) -- (6,0.75);

    \end{tikzpicture}
  }
\end{equation*}
An equivalent definition exists when the semantics is
taken to be in $\Span$, which is the following:
\begin{definition}
  Let $\NetSem{M}\in\PetriSpan$ be a guarded net with side effects. We define 
  the \emph{internalization of $\NetSem{M}$}, denoted 
  with $\Grothendieck{\Fun{M}}$, as the following category:
  \begin{itemize}
    \item The objects of $\GrothendieckS{M}$ are pairs 
    $(X, x)$, where $X$ is an object of $\Free{M}$ and $x$ is an 
    element of $\Fun{M}X$. Concisely:
    \begin{equation*}
      \Obj{\GrothendieckS{M}} := 
        \Suchthat{(X,x)}{(X \in \Obj{\Free{M}}) \wedge (x \in \Fun{M}X)}.
    \end{equation*}
    \item A morphism from $(X,x)$ to 
    $(Y,y)$ in $\Grothendieck{\Fun{M}}$ is a pair $(f,s)$ where $f\colon X \to Y$ in $\Free{M}$ and 
    $s\in S_{\Fun{M}f}$ in the apex of the corresponding span connects $x$ to $y$. Concisely:
    \begin{align*}
      &\Hom{\GrothendieckS{M}}{(X,x)}{(Y,y)} :=\\
       &\qquad :=\Suchthat{(f,s)}{(f \in \Hom{\Free{M}}{X}{Y}) \wedge (s\in S_{\Fun{M}f})\wedge(\Fun{M}f(s) = (x,y))}.
    \end{align*}
  \end{itemize}
\end{definition}
The intuition in the span case is exactly as for partial functions,
and we don't deem it useful to draw the same picture again.
Looking at the example, though, a couple of things 
become clear. The first is that to justify the
idea of the Grothendieck construction turning an assignment 
of semantics into a net we have to prove that the resulting 
category is symmetric strict monoidal and free. The second 
is that the net thus built is fibered over the base net, and there
should be an opposite construction sending $\GrothendieckS{M}$
to $M$. Both of these 
claims are true, as we now prove:
\begin{lemma}
  \label{lem: Grothendieck is a SSMC}
  In the case of both $\SetS$ and $\Span$, the category $\GrothendieckS{M}$
  has a strict symmetric monoidal structure.
\end{lemma}
\begin{theorem}
  \label{thm: Grothendieck is FSSMC}
  In both the case of $\SetS$ and of $\Span$ the strict symmetric monoidal category $\GrothendieckS{M}$ is free.
\end{theorem}
\begin{counterexample}[Relations]
  \label{counterex: single-valued is necessary}
  \cref{thm: Grothendieck is FSSMC} does not hold -- 
  the Grothendieck construction does not yield a free symmetric 
  strict monoidal category -- if we replace $\SetS$ or $\Span$ with $\Rel$. To see this, 
  consider $\Grothendieck{\Fun{M}}$ in the case that 
  $\Fun{M}\colon\Free{M} \to \Rel$.  Let $M$ be the Petri net consisting of 
  three places $X, Y, Z$ and two transitions $f\colon X \to Y$ and $g\colon Y \to Z$.
  Let $\Fun{M}$ send $X$ to $\{x\}$, $Y$ to $\{y_1, y_2\}$,  and 
  $Z$ to $\{z\}$. On morphisms, let $\Fun{M}$ send $f$ to the maximal 
  relation on $\{x\} \times \{y_1, y_2\}$ and $g$ to the maximal relation 
  on $\{y_1, y_2\} \times \{z\}$. Then we have the following four 
  generating morphisms in $\GrothendieckS{M}$:
  \begin{gather*}
   f_1\colon (X,x) \to (Y,y_1) 
    \qquad f_2\colon (X,x) \to (Y,y_2) \\
     g_1\colon (Y,y_1) \to (Z, z)
    \qquad g_2\colon (Y,y_2) \to (Z,z)
  \end{gather*}
  There is an equality $f_1\cp g_1=f_2\cp g_2$ as morphisms 
  $(X,x) \to (Z,z)$ in $\Grothendieck{\Fun{M}}$, proving 
  $\Grothendieck{\Fun{M}}$ is not free. 

  The reason that \cref{thm: Grothendieck is FSSMC} holds in the span case is that 
  spans \emph{keep track of different paths between elements, whereas relations
  do not.} To see this, consider the span composition:
  
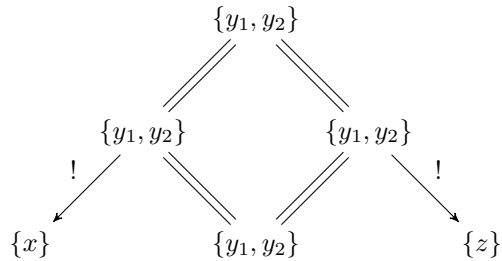
\begin{center}
    \begin{tikzpicture}[node distance=1.3cm,>=stealth',bend angle=45,auto]
      \node (left) at (0,0) {$\{x\}$};
      \node (center) at (3,0) {$\{y_1,y_2\}$};
      \node (right) at (6,0) {$\{z\}$};
      \node(leftTip) at (1.5,1.5) {$\{y_1,y_2\}$};
      \node(rightTip) at (4.5,1.5) {$\{y_1,y_2\}$};
      \node(topTip) at (3,3) {$\{y_1,y_2\}$};
      \draw[->] (leftTip) to node [midway,above left] {$!$} (left);
      \draw[transform canvas={yshift=1pt, xshift=1pt}] (leftTip) to (center);
      \draw[transform canvas={yshift=-1pt, xshift=-1pt}] (leftTip) to (center);
      \draw[transform canvas={yshift=1pt, xshift=-1pt}] (rightTip) to (center);
      \draw[transform canvas={yshift=-1pt, xshift=1pt}] (rightTip) to (center);
      \draw[->] (rightTip) to node [midway,above right] {$!$} (right);
      \draw[transform canvas={yshift=1pt, xshift=-1pt}] (topTip) to (leftTip);
      \draw[transform canvas={yshift=-1pt, xshift=1pt}] (topTip) to (leftTip);
      \draw[transform canvas={yshift=1pt, xshift=1pt}] (topTip) to (rightTip);
      \draw[transform canvas={yshift=-1pt, xshift=-1pt}] (topTip) to (rightTip);
    \end{tikzpicture}
    \end{center} 
    It is clear that in this composition the two paths from $x$ to $z$
    are considered as separated in the $\Span$ case, and witnessed by $y_1,y_2$ respectively, 
    while in the case of $\Rel$ they would have been 
    conflated to one.  The result is that these paths correspond 
    to the same morphism in the relational case of $\GrothendieckS{M}$, 
    introducing new equations and breaking freeness, while 
    they stay separated in the span case.
\end{counterexample}
\begin{lemma}
  \label{lem: strict monoidal functor pi SetS}
  In the case of both $\SetS$ and $\Span$, there is a strict monoidal functor 
  $\pi_M\colon \GrothendieckS{M} \to \Free{M}$ 
   sending $(X,x)$ to $X$ and $f\colon(X,x) \to (Y,y)$ to $f\colon X \to Y$
   (resp.\ $(f,s)\colon(X,x) \to (Y,y)$ to $f\colon X \to Y$).
\end{lemma}
\begin{remark}\label{rem: pi not full}
 In general, $\pi_M$ is not an opfibration. This is because our 
 target categories $\SetS$ and $\Span$ allow for partial functions. 
 Indeed, if $f\colon X\to Y$ in $M$ is sent by $\Fun{M}$ to a partial function 
 that is not defined on $x\in \Fun{M}X$, 
 then there is no coCartesian lift emanating from $(X,x)$ for the morphism $f$.
\end{remark}
We conclude this section by proving that the reachability semantics
of a guarded net coincides with the reachability semantics of its internalization.
\begin{theorem}
  \label{thm: reachability is preserved}
  Let $\NetSem{N}$ be a guarded Petri net (with side effects).
  $(Y,y)$ is reachable from $(X,x)$ if and only if $(Y,y)$ is reachable
  from $(X,x)$ in the net $\UnFree{\GrothendieckS{N}}$.
\end{theorem}
\section{Properties of internalizations}
The Grothendieck construction provides a way to internalize 
partial function and span semantics to nets. As such, it acts 
on objects of the categories $\PetriSetS$ and $\PetriSpan$, 
respectively. It is thus worth asking what happens to morphisms 
in these categories. The answer is, luckily, easy to find:
\begin{lemma}
  \label{lem: lift F}
  Let $F\colon \NetSem{M} \to \NetSem{N}$ be a morphism in $\PetriSetS$ (resp. in $\PetriSpan$). Then 
  it lifts to strict monoidal functor
  $\LiftSetS{F}\colon \GrothendieckS{M} \to \GrothendieckS{N}$ 
  (resp.\ 
  $\LiftSpan{F}\colon \GrothendieckS{M} \to \GrothendieckS{N}$), 
  such that the following diagram on the left (resp.\ on the right) commutes:
  \begin{equation*}
    \begin{tikzpicture}[node distance=1.3cm,>=stealth',bend angle=45,auto]
      \node (1) at (0,1.5) {$\GrothendieckS{M}$};
      \node (2) at (0,0) {$\Free{M}$};
      \node (3) at (3,1.5) {$\GrothendieckS{N}$};
      \node (4) at (3,0) {$\Free{N}$};
      \node (5) at (1.5, -1.5) {$\SetS$};
      \draw[->] (1)--(2) node [midway,left] {$\pi_M$};
      \draw[->] (3)--(4) node [midway,right] {$\pi_N$};
      \draw[->] (2)--(4) node [midway,above] {$F$};
      \draw[->, dashed] (1)--(3) node [midway,above] {$\LiftSetS{F}$};
      \draw[->] (2)--(5) node [midway,left] {$\Fun{M}$};
      \draw[->] (4)--(5) node [midway,right] {$\Fun{N}$};
    \end{tikzpicture}
    \hspace{6em}
    \begin{tikzpicture}[node distance=1.3cm,>=stealth',bend angle=45,auto]
      \node (1) at (0,1.5) {$\GrothendieckS{M}$};
      \node (2) at (0,0) {$\Free{M}$};
      \node (3) at (3,1.5) {$\GrothendieckS{N}$};
      \node (4) at (3,0) {$\Free{N}$};
      \node (5) at (1.5, -1.5) {$\Span$};
      \draw[->] (1)--(2) node [midway,left] {$\pi_M$};
      \draw[->] (3)--(4) node [midway,right] {$\pi_N$};
      \draw[->] (2)--(4) node [midway,above] {$F$};
      \draw[->, dashed] (1)--(3) node [midway,above] {$\LiftSpan{F}$};
      \draw[->] (2)--(5) node [midway,left] {$\Fun{M}$};
      \draw[->] (4)--(5) node [midway,right] {$\Fun{N}$};
    \end{tikzpicture}
  \end{equation*}
\end{lemma}
\begin{notation}
  The notation for the liftings in \cref{lem: lift F}
  is easy to remember: The arrow over $\LiftSetS{F}$ looks
  like a stylized function, while the hat over $\LiftSpan{F}$ 
  looks like a sylized span.
\end{notation}
The lifting of \cref{lem: lift F} is quite well-behaved. 
First of all, it is worth stressing how it preserves some 
relevant categorical properties:
\begin{lemma}
  \label{lem: hatF SetS is full and faithful}
  For any map $F\colon(M,\Fun{M})\to(N,\Fun{N})$ in $\PetriSetS$ (respectively in $\PetriSpan$), 
  the functor $F$ is faithful if and only if $\LiftSetS{F}$ is faithful (resp.\ $\LiftSpan{F}$ is faithful). 
  If $F$ is full, then so is $\LiftSetS{F}$ (resp.\ $\LiftSpan{F}$).  
\end{lemma}
Having ascertained that ``basic'' categorical properties are preserved,
it is worth asking what happens to particular classes of functors in 
$\PetriSetS$ and $\PetriSpan$, respectively. 

Following~\cite{Genovese2019}, there are three relevant kinds of morphisms in 
a category of Petri nets with semantics.
On one hand there are transition-preserving 
functors, which represent morphisms of free monoidal categories 
arising purely from the topological structure of their underlying net.
On the other there are functors representing glueings 
of nets, which are themselves divided into 
synchronizations (defined in terms of addition and erasing
of generators, that is, double pushouts) and identifications 
(defined in terms of pushouts). Let us 
investigate which ones of these properties are preserved.
\begin{definition}
  A strict symmetric monoidal functor $F$ between
  FSSMCs $\CategoryC$, $\CategoryD$ is said to be \emph{transition-preserving}
  when each generating morphism $f$ of $\CategoryC$ is
  mapped to $\sigma\cp g\cp \sigma'$ for some generating 
  morphism  $g$  of $\CategoryD$ and symmetries $\sigma, \sigma'$.
\end{definition}
\begin{lemma}
  \label{lem: hat F SetS is transition-preserving}
    If $F$ is transition-preserving, so are $\LiftSetS{F}$ and $\LiftSpan{F}$.
\end{lemma}
\begin{lemma}
  \label{lem: hat F SetS is injective on objects}
  If $F$ is injective on objects, so are $\LiftSetS{F}$ and $\LiftSpan{F}$.
\end{lemma}
Identifications are also preserved. The ultimate reason for this
is that identifications are defined purely in terms of coequalizers of transition-preserving
functors, which are preserved by the Grothendieck construction.
\begin{definition}
  A Petri net $\NetSem{N}$ is said to be an 
  \emph{identification of $\NetSem{M}$} if
  there is a morphism $F: \NetSem{M} \to \NetSem{N}$ such that:
  \begin{itemize}[itemsep=0pt]
    \item There is a Petri net $O$, and a pair of 
    \emph{transition-preserving} functors 
    $l,r: \Free{O} \to \Free{M}$;
    \item $l\cp\Fun{M} = r\cp\Fun{M}$; and
    \item $F$ is the coequalizer of $l$ and $r$.
  \end{itemize}
\end{definition}
\begin{lemma}\label{lem: hat F preserves identifications}
  If $\NetSem{N}$ is an identification of $\NetSem{M}$ via
  $F$ and witnesses $O, l, r$, then $\GrothendieckS{N}$ is an 
  identification of $\GrothendieckS{M}$ via $\LiftSetS{F}$ and 
  witnesses $\UnFree{\Grothendieck{(l \cp \Fun{M})}}, \LiftSetS{l}, \LiftSetS{r}$.
  The span case is analogous.
\end{lemma}
Preservation of identifications also entails that both
additions and erasings of generators for a net are preserved by
internalizations.
\begin{definition}
  A net $\NetSem{M}$ is an
  \emph{addition of generating morphisms to $\NetSem{K}$
  via $W, w$} if: 
  \begin{itemize}
    \item  There is a net $W$ together with a strict monoidal functor 
    $w: \Free{W} \to \Free{N}$ which sends 
    generating objects to generating objects, 
    is injective on objects and faithful;
    \item $\Free{M}$ is the pushout of 
    $\Free{\overline{W}} \hookrightarrow \Free{W} \xrightarrow{w} \Free{K}$
    and $\Free{\overline{W}} \hookrightarrow \Free{W}$, 
    where $\overline{W}$ denotes the net with the same places of $W$ and no transitions; and
    \item $\Fun{M}$ arises from
    the universal property of the pushout.
  \end{itemize}
\end{definition}
\begin{lemma}
  \label{lem: hat F preserves addition of generating morphisms}
  Let $\NetSem{M}$  be an addition of generating
  morphisms to $\NetSem{K}$ via $W, w$. 
  Then $\GrothendieckS{M}$ is an addition of 
  generating morphisms to $\GrothendieckS{K}$ via 
  $\UnFree{\Grothendieck{(w \cp \Fun{K})}},
  \LiftSetS{w}$. The span case is analogous.
\end{lemma}
\begin{definition}
  Let $(N_w, \Fun{(N_w)})$ be a subnet of $\NetSem{N}$. An \emph{erasing of 
  generators of $\NetSem{N}$ via $N_w$} is a net 
  $\NetSem{K}$  such that:
  \begin{itemize}[itemsep=0pt]
    \item $\NetSem{K}$ is a subnet of $\NetSem{N}$;
    \item $\NetSem{\overline{N_w}}$, where $\overline{N_w}$ 
    denotes the net with the same places of $N_w$ and no transitions, is a subnet of $\NetSem{K}$; and
    \item $\Free{N}$ is the pushout of $\Free{\overline{N_w}} \hookrightarrow \Free{N_w}$;
    and $\Free{\overline{N_w}} \hookrightarrow \Free{K}$.
  \end{itemize}
\end{definition}
\begin{lemma}
  \label{lem: hat F preserves erasing of generating morphisms}
  Let $\NetSem{K}$ be an erasing of generating
  morphisms from $\NetSem{N}$ via a subnet $N_w$.
  Then $\GrothendieckS{K}$ is an erasing 
  of generators from $\GrothendieckS{N}$ via
  $\Grothendieck sub_{N_w} \cp \Fun{N}$. The span case is 
  analogous.
\end{lemma}
Surprisingly, even if erasing and addition of generators
are preserved by internalizations, synchronizations are not.
Indeed, following~\cite{Genovese2019}, $\NetSem{M}$ is a synchronization 
of $\NetSem{N}$ via $W, w$ when $\Free{M}$ is defined to be the 
result of applying a double
pushout rewrite rule $\Free{N_w} \xleftarrow{w'} \Free{W} 
\xhookleftarrow{in_W} \Free{\overline{W}} 
\xhookrightarrow{in_W} \Free{W}$ to $\Free{N}$, 
with $w$ factorizing through $w'$.
In internalizing this construction the pushouts are 
preserved, but the rewrite rule is not!  This becomes evident by lifting
the definition of synchronization altogether, where
in the following diagram we are sticking to the
notation developed in~\cite{Genovese2019}:
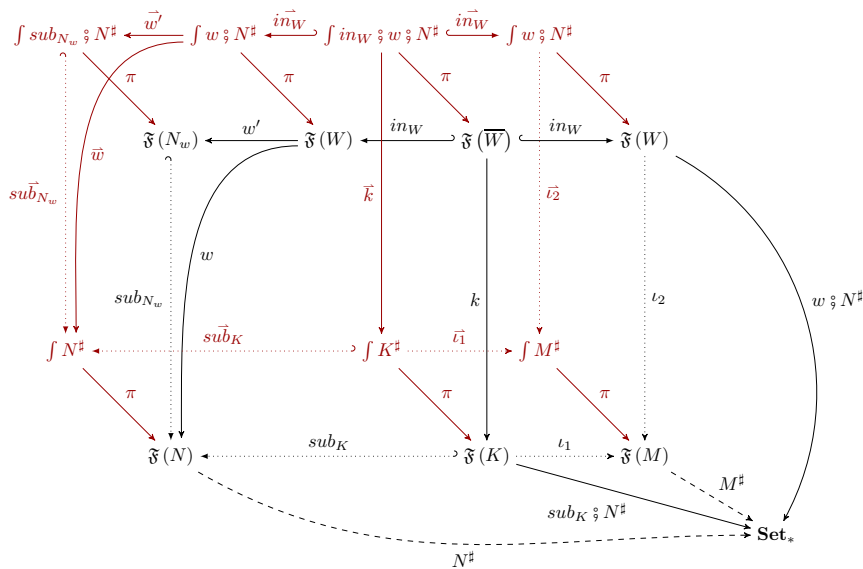
\begin{equation*}
  \scalebox{0.7}{
    \begin{tikzpicture}[node distance=2cm,>=stealth',bend angle=45,auto]
      \node (in1) at (0,6) {$\Free{N_w}$};
      \node (in2) at (3,6) {$\Free{W}$};
      \node (in3) at (6,6) {$\Free{\overline{W}}$};
      \node (in4) at (9,6) {$\Free{W}$};
      \node (in5) at (0,0) {$\Free{N}$};
      \node (in6) at (6,0) {$\Free{K}$};
      \node (in7) at (9,0) {$\Free{M}$};
  
      \node (sem) at (11.5, -1.5) {$\SetS$};
  
      \node[red!60!black] (out1) at (-2,8) {$\Grothendieck{sub_{N_w} \cp \Fun{N}}$};
      \node[red!60!black] (out2) at (1,8) {$\Grothendieck{w \cp \Fun{N}}$};
      \node[red!60!black] (out3) at (4,8) {$\Grothendieck{in_W \cp w \cp \Fun{N}}$};
      \node[red!60!black] (out4) at (7,8) {$\Grothendieck{w \cp \Fun{N}}$};
      \node[red!60!black] (out5) at (-2,2) {$\Grothendieck{\Fun{N}}$};
      \node[red!60!black] (out6) at (4,2) {$\Grothendieck{\Fun{K}}$};
      \node[red!60!black] (out7) at (7,2) {$\Grothendieck{\Fun{M}}$};
  
      \draw[left hook-latex] (in3) to node[swap] {$in_W$} (in2);
      \draw[->] (in2) to node[swap] {$w'$} (in1);
      \draw[right hook-latex] (in3) to node {$in_W$} (in4);
      \draw[left hook-latex, dotted] (in1) to node[swap] {$sub_{N_w}$} (in5);
      \draw[->, out=180, in=90] ([yshift=2mm]in2.south west) to node {$w$} ([xshift=2mm, yshift=3.5mm]in5);
      \draw[->] (in3) to node[left] {$k$} (in6);
      \draw[left hook-latex, dotted] (in6) to node[swap] {$sub_K$} (in5);
      \draw[->, dotted] (in6) to node {$\iota_1$} (in7);
      \draw[->, dotted] (in4) to node {$\iota_2$} (in7);
  
      \draw[left hook-latex, red!60!black] (out3) to node[swap] {$\LiftSetS{in_W}$} (out2);
      \draw[->, red!60!black] (out2) to node[swap] {$\LiftSetS{w'}$} (out1);
      \draw[right hook-latex, red!60!black] (out3) to node {$\LiftSetS{in_W}$} (out4);
      \draw[left hook-latex, dotted, red!60!black] (out1) to node[swap] {$\LiftSetS{sub_{N_w}}$} (out5);
      \draw[->, red!60!black, out=180, in=90] ([yshift=2mm]out2.south west) to node {$\LiftSetS{w}$} ([xshift=2mm, yshift=3.5mm]out5);
      \draw[->, red!60!black] (out3) to node[left] {$\LiftSetS{k}$} (out6);
      \draw[left hook-latex, dotted, red!60!black] (out6) to node[swap] {$\LiftSetS{sub_K}$} (out5);
      \draw[->, dotted, red!60!black] (out6) to node {$\LiftSetS{\iota_1}$} (out7);
      \draw[->, dotted, red!60!black] (out4) to node {$\LiftSetS{\iota_2}$} (out7);
  
      \draw[->, red!60!black] (out1) to node {$\pi$} (in1);
      \draw[->, red!60!black] (out2) to node {$\pi$} (in2);
      \draw[->, red!60!black] (out3) to node {$\pi$} (in3);
      \draw[->, red!60!black] (out4) to node {$\pi$} (in4);
      \draw[->, red!60!black] (out5) to node {$\pi$} (in5);
      \draw[->, red!60!black] (out6) to node {$\pi$} (in6);
      \draw[->, red!60!black] (out7) to node {$\pi$} (in7);
  
      \draw[->, bend left] (in4) to node {$w \cp \Fun{N}$} (sem);
      \draw[->, dashed,  out=-30, in= 180] (in5) to node[swap] {$\Fun{N}$} (sem);
      \draw[->] (in6) to node[swap] {$sub_K \cp \Fun{N}$} (sem);
      \draw[->, dashed] (in7) to node {$\Fun{M}$} (sem);
  
    \end{tikzpicture}
  }
\end{equation*}
The black arrows are just the definition of synchronization.
The dotted arrows denote the pushout arrows, while the
dashed arrows arise from the universal property of the pushout.
The maroon arrows and objects represent the Grothendieck construction
and the lifting of the functors obtained from \cref{lem: lift F},
while the $\pi$s stand for the functors obtained in 
\cref{lem: strict monoidal functor pi SetS}, where we omitted 
subscripts to avoid clutter. 

As one can see the pushout squares are both preserved,
but $\GrothendieckS{M}$ is not a synchronization of
$\GrothendieckS{N}$ via $\Grothendieck{in_W \cp w \cp \Fun{N}}$
since
\begin{equation*}
  \Grothendieck{sub_{N_w} \cp \Fun{N}} \neq \left(\GrothendieckS{N}\right)_{\Grothendieck{w \cp \Fun{N}}}
\end{equation*}
In other words, $\Grothendieck{sub_{N_w} \cp \Fun{N}}$
is too big of a subcategory of $\GrothendieckS{N}$ to make
$\GrothendieckS{M}$ into a synchronization. An analogous
observation holds for spans.

\begin{counterexample}[Synchronizations not preserved]
  We provide a practical counterexample of why
  synchronizations are not preserved by internalizations.
  Consider the following nets, where we 
  are borrowing the graphical notation developed in~\cite{Genovese2019}, 
  decorating net elements with their images in $\SetS$.
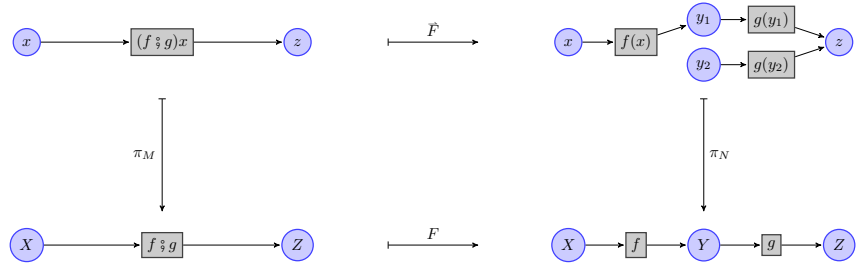
  \begin{equation*}
    \scalebox{0.6}{
      \begin{tikzpicture}[node distance=1.3cm,>=stealth',bend angle=45,baseline=2cm]
        \node [place, tokens=0,] (1aBL) at (0,-1.5){$X$};
        \node [transition] (2aBL) at (3,-1.5)     {$f \cp g$}
          edge [pre] (1aBL);
        \node [place,tokens=0] (3aBL) at (6,-1.5)  {$Z$}
          edge [pre] (2aBL);
        \draw[|->,] (3,1.75) to node[midway, left] (piM) {$\pi_M$} (3,-0.75);
        \node [place, tokens=0,] (1aTL) at (0,3){$x$};
        \node [transition] (2aTL) at (3,3)     {$(f \cp g)x$}
          edge [pre] (1aTL);
        \node [place,tokens=0] (3aTL) at (6,3)  {$z$}
          edge [pre] (2aTL);
        \draw[|->,] (8,-1.5) to node[midway, above] (F) {${F}$} (10,-1.5);
        \draw[|->,] (8,3) to node[midway, above] (LiftSetSF) {$\LiftSetS{F}$} (10,3);
        \node [place, tokens=0,] (1aBR) at (12,-1.5){$X$};
        \node [transition] (2aBR) at (13.5,-1.5)     {$f$}
          edge [pre] (1aBR);
        \node [place, tokens=0] (3aBR) at (15,-1.5)     {$Y$}
          edge [pre] (2aBR);
        \node [transition] (4aBR) at (16.5,-1.5)     {$g$}
          edge [pre] (3aBR);
        \node [place,tokens=0] (5aBR) at (18,-1.5)  {$Z$}
          edge [pre] (4aBR);
        \draw[|->,] (15,1.75) to node[midway, right] (piN) {$\pi_N$} (15,-0.75);
        \node [place, tokens=0,] (1aTR) at (12,3){$x$};
        \node [transition] (2aTR) at (13.5,3)     {$f(x)$}
          edge [pre] (1aTR);
        \node [place, tokens=0,] (3aTR) at (15,3.5){$y_1$}
          edge [pre] (2aTR);
        \node [place, tokens=0,] (3bTR) at (15,2.5){$y_2$};
        \node [transition] (4aTR) at (16.5,3.5)     {$g(y_1)$}
          edge [pre] (3aTR);
        \node [transition] (4bTR) at (16.5,2.5)     {$g(y_2)$}
          edge [pre] (3bTR);
        \node [place,tokens=0] (5aTR) at (18,3)  {$z$}
          edge [pre] (4aTR)
          edge [pre] (4bTR);
      \end{tikzpicture}
    }
  \end{equation*}
  At the base level we have two nets, $M$ on the left 
  and $N$ on the right.  Eliding the functor $\Fun{N}$, the places of $N$ are mapped to sets:
  \begin{equation*}
    X \coloneqq \{x\} \qquad Y \coloneqq \{y_1, y_2\} \qquad Z \coloneqq \{z\}
  \end{equation*}
  While transitions are mapped to partial functions $f\colon X \to Y$ 
  and $g\colon Y \to Z$, defined as follows:
  \begin{equation*}
    f(x) = y_1 \qquad g(y_1) = g(y_2) = z
  \end{equation*}
  $M$ is clearly a synchronization of $N$ via $F$:
  The generators $f$ and $g$ have been erased and
  a generator corresponding to $f\cp g$ has been added.
  Taking the Grothendieck construction on $M$ and $N$
  (top left and top right in the figure, respectively),
  we see how the erasing of generators is problematic:
  The morphism $g$ in $N$ branches into $g(y_1)$
  and $g(y_2)$ in $\GrothendieckS{N}$, 
  of which only $g(y_1)$ forms a path 
  with $f(x)$. In lifting the synchronization $M$ to 
  $\GrothendieckS{M}$, we would expect 
  $g(y_1)$ and $f(x)$ to be erased and conflated
  into $f(x)\cp g(y_1)$, whereas $g(y_2)$ stays.
  But this is not the case, since in $M$ the generator
  $g$ has already been erased ``before being allowed to branch'', 
  taking $g(y_2)$ with it when we take $\GrothendieckS{M}$!

  As we said before, this ultimately depends on the
  fact that the internalization of
  the subnet provided by the synchronization
  witness contains too many morphisms,
  and ends up erasing more generators than
  we would like it to.
\end{counterexample}
\section{Internalization as a functor}
In this final section, we put together some of the
properties we have proved so far about internalizations,
and prove that internalization is a functor. 
The intuitive argument behind the results that are about
to follow is this: If $\GrothendieckS{N}$ internalizes the semantics
of $\NetSem{N}$, in either the case $\Fun{N}\colon N\to\SetS$ or 
$\Fun{N}\colon N\to\Span$, then $\GrothendieckS{N}$ should be considered as
``just a net", that is, an object of $\PetriTerm$; see \cref{ex.petri_one}.
 
Putting together results about lifting of functors
obtained in the previous section, we are indeed
able to prove this.
\begin{theorem}
  \label{thm: mapping to Term}
  Denote with $\TermCat$ the terminal category, together with the
   trivial symmetric monoidal structure on it. There is a faithful, 
   strong monoidal functor $\EmbSetS\colon\PetriSetS \to \PetriTerm$
   defined as follows:
  \begin{itemize}
    \item On objects, it sends $\NetSem{M}$ to 
    $\NetSemTerm{\UnFree{\GrothendieckS{M}}}$.
    \item On morphisms, we send the functor
      $F\colon \NetSem{M} \to \NetSem{N}$ to the 
      functor\footnote{
        To be absolutely precise, we are referring to
        the functor $\Free{\UnFree{\GrothendieckS{M}}} \simeq 
        \GrothendieckS{M} \xrightarrow{\hat{F}} \GrothendieckS{N} 
        \simeq \Free{\UnFree{\GrothendieckS{N}}}$.}
      \begin{equation*}
        \hat{F}\colon \NetSemTerm{\UnFree{\GrothendieckS{M}}} 
          \to  \NetSemTerm{\UnFree{\GrothendieckS{N}}}
      \end{equation*}
  \end{itemize}
  Similarly, there is a faithful, 
  strong monoidal functor
  $\EmbSpan\colon\PetriSpan \to \PetriTerm$
  defined as follows:
  \begin{itemize}
    \item On objects, it sends $\NetSem{M}$ to 
    $\NetSemTerm{\UnFree{\GrothendieckS{M}}}$.
    \item On morphisms, we send the functor
      $F\colon\NetSem{M} \to \NetSem{N}$ to the 
      functor\footnote{
        To be absolutely precise, we are referring to
        the functor $\Free{\UnFree{\GrothendieckS{M}}} \simeq 
        \GrothendieckS{M} \xrightarrow{\hat{F}} \GrothendieckS{N} 
        \simeq \Free{\UnFree{\GrothendieckS{N}}}$.}
      \begin{equation*}
        \hat{F}\colon \NetSemTerm{\UnFree{\GrothendieckS{M}}} 
          \to  \NetSemTerm{\UnFree{\GrothendieckS{N}}}
      \end{equation*}
  \end{itemize}
\end{theorem}
Finally, it is worth nothing that for each
choice of semantics $\Semantics$ there is 
another obvious functor from $\PetriS{\Semantics}$ 
to $\PetriTerm$, which just forgets the semantics
altogether. It is worth asking how this functor
and the ones provided in \cref{thm: mapping to Term}
are related.
\begin{proposition}
  \label{lem: natural transformations from emb to for}
  Denote with 
  \begin{equation*}
    \ForSetS\colon\PetriSetS \to \PetriTerm \qquad \ForSpan\colon
    \PetriSpan \to \PetriTerm
  \end{equation*}
  the ``forgetful'' functors defined 
  by sending each Petri net $\NetSem{M}$
  to $\NetSemTerm{M}$.
  Then there are natural transformations:
  \begin{equation*}
    \begin{tikzpicture}[node distance=1.3cm,>=stealth',bend angle=45,auto]
      \node (1) at (0,0) {$\PetriSetS$};
      \node (2) at (3,0) {$\PetriTerm$};
      \node[rotate=-90] (3) at (1.5,0) {$\Rightarrow$};
      \node(4) at (1.8,0) {$\pi$};
      \draw[->, bend right] (1) to node [midway,below] {$\ForSetS$} (2);
      \draw[->, bend left] (1) to node [midway,above] {$\EmbSetS$} (2);
    \end{tikzpicture}
    \qquad
    \begin{tikzpicture}[node distance=1.3cm,>=stealth',bend angle=45,auto]
      \node (1) at (0,0) {$\PetriSpan$};
      \node (2) at (3,0) {$\PetriTerm$};
      \node[rotate=-90] (3) at (1.5,0) {$\Rightarrow$};
      \node(4) at (1.8,0) {$\pi$};
      \draw[->, bend right] (1) to node [midway,below] {$\ForSpan$} (2);
      \draw[->, bend left] (1) to node [midway,above] {$\EmbSpan$} (2);
    \end{tikzpicture}
  \end{equation*}
\end{proposition}
\section{Conclusion and future work}
In this work, we described guarded Petri nets
as Petri nets endowed with a functorial semantics. We provided
two different styles of semantics: a deterministic one,
realized using the category of sets and partial functions, and
a non-deterministic one that allows for side effects, realized using the category of
partial functions and that of spans.

We moreover showed how, using the Grothendieck construction,
the guards can be internalized, obtaining a Petri net whose
reachability relation is equivalent to the one of the guarded one.
We proved that internalizations have nice properties, and the internalization
construction is functorial in the choice of the guarded net we start from.

Regarding directions of future work, a pretty straightforward thing to
do would be to figure out which semantics, other than $\SetS$ and 
$\Span$, are internalizable. That is, if $F\colon \Free{N} \to \Semantics$ is
a symmetric monoidal functor, which properties do $\Semantics$ and $F$
need to have so that $\Grothendieck{F}$ is a free symmetric strict monoidal category.
\section{Acknowledgements}
David Spivak acknowledges support from Honeywell Inc.\ as well as from AFOSR grants FA9550-17-1-0058 and FA9550-19-1-0113.
Fabrizio Genovese wants to thank his fellow team members at Statebox for useful discussion and support.
\bibliographystyle{splncs04}
\bibliography{Bibliography}
\clearpage
\appendix\label{appendix}
\section*{Appendix -- Proofs}
\begingroup
\def\thelemma{\ref{lem: Grothendieck is a SSMC}}
\begin{lemma}
  In the case of both $\SetS$ and $\Span$, the category $\GrothendieckS{M}$
  has a strict symmetric monoidal structure.
\end{lemma}
\addtocounter{lemma}{-1}
\endgroup
\begin{proof}
  We start with the case of $\SetS$. Since $\Fun{M}$ is strict monoidal, $\Fun{M}(X\otimes Y)=\Fun{M}X\times \Fun{M}Y$.
  Thus on objects, we can set $(X,x) \Tensor (Y,y) \coloneqq (X \Tensor Y, (x,y))$. On morphisms, 
  we just use the monoidal product $f \Tensor g$ from $\Free{M}$. The monoidal unit is $(\TensorUnit, *)$,
  where $I$ is the monoidal unit of $\Free{M}$ and $*$ 
  is the unique element of the monoidal unit $\{*\}$ of $\SetS$. 

  First, we have to prove that  $(\_) \Tensor (\_)$  is well-defined, 
  namely that if $f\colon (X,x) \to (X',x')$ and $g\colon (Y,y) \to (Y',y')$ then 
  $\Fun{M}(f \Tensor g)(x,y) = (x', y')$. This is obvious since 
  $\Fun{M}$ is strict, therefore:
  \begin{equation*}
    \Fun{M}(f \Tensor g)(x,y) 
      = (\Fun{M}f \times \Fun{M}g)(x,y) 
      = (\Fun{M}f(x), \Fun{M}g(y)) 
      = (x',y').
  \end{equation*}
  Bifunctoriality of $(\_) \Tensor (\_)$ follows trivially from the fact
  that it is defined as in $\Free{M}$. Moreover, for each object $(X,x)$ 
  we have:
  \begin{gather*}
    (X, x) \Tensor (\TensorUnit, *) = (X \Tensor \TensorUnit, (x,*)) = (X, (x,*)) = (X,x)\\
    (\TensorUnit, *) \Tensor (X, x) = (\TensorUnit \Tensor X, (*,x)) = (X, (*,x)) = (X,x)
  \end{gather*}
  Where the rightmost equality in each line follows from the fact 
  that we are working with the strictified version of $\SetS$. 
  This proves that unitors are identities. A similar reasoning 
  can be applied to associators, concluding the proof.

  \medskip
  \noindent
  Now we consider the case of $\Span$.
  On objects, we set again $(X,x) \Tensor (Y,y) \coloneqq (X \Tensor Y, (x,y))$. 
  On morphisms, we set $(f,s) \Tensor (g,t)\coloneqq(f\otimes g, (s,t))$, where 
  $f\otimes g$ is as in $\Free{M}$ and $(s,t)$ is the pair of span-apex elements:
  \[
  \begin{tikzcd}[row sep=0]
    \Fun{M}X&
    S\ar[l, "f_1"']\ar[r, "f_2"]&
    \Fun{M}X'\\
    (x,y)&
    (s,t)\ar[l, "{(f_1,g_1)}"']\ar[r, "{(f_2,g_2)}"]&
    (x',y')\\
    \Fun{M}Y&
    T\ar[l, "g_1"']\ar[r, "g_2"]&
    \Fun{M}Y'
  \end{tikzcd}
  \]
  The monoidal unit is $(\TensorUnit, *)$,
  where $I$ is the monoidal unit of $\Free{M}$ and $*$ 
  is the unique element of the monoidal unit $\{*\}$ of $\Span$. 

  First, we have to prove that the monoidal product  $(\_) \Tensor (\_)$  is well-defined, 
  namely that if $(f,s)\colon (X,x) \to (X',x')$ and $(g,t)\colon (Y,y) \to (Y',y')$ then 
  $\Fun{M}(f \Tensor g)(s,t)=((x,y) ,(x', y'))$. This is obvious since 
  $\Fun{M}$ is strict, therefore:
  \begin{equation*}
    \Fun{M}(f \Tensor g)(s,t) 
      = (\Fun{M}f \times \Fun{M}g)(s,t) 
      = (\Fun{M}f(s), \Fun{M}g(t)) 
      = ((x,y),(x',y')).
  \end{equation*}
  The remainder of the proof is as in the previous case.
\end{proof}
\begingroup
\def\thetheorem{\ref{thm: Grothendieck is FSSMC}}
\begin{theorem}
  In both the case of $\SetS$ and of $\Span$ the strict symmetric monoidal category $\GrothendieckS{M}$ is free.
\end{theorem}
\addtocounter{theorem}{-1}
\endgroup
\begin{proof}
We start with the case of $\SetS$.  Consider the free symmetric strict monoidal category $\hat{M}$ generated as follows:
  \begin{itemize}
    \item Object generators are pairs $(X,x)$, with $X$ object generator in 
    $\Free{M}$ and $x \in \Fun{M}X$;
    \item A morphism generator $(X,x) \to (Y,y)$ is a morphism generator $f\colon X \to Y$ 
    of $\Free{M}$ such that $\Fun{M}f(x) = y$.
  \end{itemize}
  We want to prove that $\hat{M}$ and $\GrothendieckS{M}$ are isomorphic.

  First, let $(X,x)$ be an object in $\GrothendieckS{M}$. Then $X$ is an object of $\Free{M}$,
  which is free, and hence we have $X = X_1 \Tensor \dots \Tensor X_n$ for 
  generating objects $X_1 \dots X_n$ in $\Free{M}$. By definition, we have 
  $x \in \Fun{M}X$. Being $\Fun{M}$ strict, this means:
  \begin{align*}
    x \in \Fun{M}X
      &\Leftrightarrow x \in \Fun{M}(X_1 \Tensor \dots \Tensor X_n)\\
      &\Leftrightarrow x \in \Fun{M}X_1 \times \dots \times \Fun{M}X_n\\
      &\Leftrightarrow \exists!(x_1 \in \Fun{M}X_1), \dots,\exists! (x_n \in \Fun{M}X_n).
       ( x = (x_1, \dots, x_n))
  \end{align*}
  Hence $(X,x)$ = $(X_1, x_1) \Tensor \dots \Tensor (X_n, x_n)$, and the object generators of 
  $\GrothendieckS{M}$ are the pairs $(X, x)$ with $X$ object generator in 
  $\Free{M}$ and $x \in \Fun{M}X$. 
  This means that there is a bijection on objects of $\hat{M}$ and $\GrothendieckS{M}$:
  \begin{equation*}
    (X_1, x_1) \Tensor \dots \Tensor (X_n, x_n) \mapsto (X_1 \Tensor \dots \Tensor X_n, (x_1, \dots, x_n))
  \end{equation*}
  We can then define a symmetric monoidal functor 
  $T\colon\hat{M} \to \GrothendieckS{M}$ that is 
  bijective on objects, and sends a generating morphism $f\colon (X,x) \to (Y,y)$ of
  $\hat{M}$ to the morphism $f\colon (X,x) \to (Y,y)$ in $\GrothendieckS{M}$.

  We want to prove that $T$ is full and faithful.
  Faithfulness is obvious; given $f_1, f_2\colon (X,x) \to (Y,y)$ in $\hat{M}$, 
  if $T(f_1)=T(f_2)$ then in particular $f_1=f_2$ in $\Free{M}$. It follows 
  from the fact that $\Fun{M}$ is at most single-valued -- i.e.\ $\Fun{M}g(x)=x'$ and $\Fun{M}g(x)=x''$ 
  imply $x'=x''$ -- that $f_1=f_2$ also in $\hat{M}$.

 For fullness, take a morphism $f\colon (X,x) \to (Z,z)$ in 
  $\GrothendieckS{M}$, and notice the following:
  \begin{itemize}
    \item If $f\colon (X,x) \to (Z,z)$ is equal to $f_1 \cp f_2$, where $f_1\colon X \to Y$ 
    and $f_2: Y \to Z$, then we have:
    \begin{equation*}
      z =\Fun{M}f (x)
        = \Fun{M}(f_1 \cp f_2)  (x)
        = (\Fun{M}f_1 \cp \Fun{M}f_2) (x)
    \end{equation*}
    So there is a $y \in \Fun{M}Y$ such that $\Fun{M}f_1(x) = y$ and 
    $\Fun{M}f_2(y) = z$. This means that $f_1\colon (X,x) \to (Y,y)$ and 
    $f_2\colon (Y,y) \to (Z,z)$ are morphisms in $\GrothendieckS{M}$;
    \item If $f\colon (X_1 \Tensor X_2,(x_1, x_2)) \to (Y_1 \Tensor Y_2, (y_1, y_2))$ 
    is equal to $f_1 \Tensor f_2$, where $f_1\colon X_1 \to Y_1$ and $f_2: X_2 \to Y_2$, 
    then we have:
    \begin{equation*}
      (y_1,y_2)= \Fun{M}f (x_1, x_2)
        = \Fun{M}(f_1 \Tensor f_2) (x_1, x_2)
        = (\Fun{M}f_1 \times \Fun{M}f_2)(x_1, x_2)
    \end{equation*}
    This means that $\Fun{M}f_1(x_1) = y_1$ and 
    $\Fun{M}f_2(x_2) = y_2$, and hence that
    $f_1\colon (X_1,x_1) \to (Y_1,y_1)$ and 
    $f_2\colon (X_2,x_2) \to (Y_2,y_2)$ are morphisms in $\GrothendieckS{M}$.
  \end{itemize}
  By definition, since $\Free{M}$ is free, any morphism $f\colon X\to Z$ can 
  be decomposed into a composition of monoidal products of morphism generators, 
  symmetries and identities. The points above prove that $f\colon (X,x)\to (Z,z)$ can be 
  decomposed in the same way, and hence is in the image of $T$; thus it is full.

  Our correspondence is bijective on objects 
  and fully faithful, proving that $\hat{M}$ and $\GrothendieckS{M}$
  are isomorphic as categories. Since $\hat{M}$ is free so is $\GrothendieckS{M}$, 
  completing the proof.

  \medskip
  \noindent
  We now consider the case of $\Span$. The structure
  of the proof is similar.
  Consider the free symmetric strict monoidal category $\hat{M}$ generated as follows:
  \begin{itemize}
    \item Object generators are pairs $(X,x)$, with $X$ object generator in 
    $\Free{M}$ and $x \in \Fun{M}X$;
    \item For each morphism generator $f\colon X \to Y$ and $s \in S$ such that 
    $\Fun{M}f(s) = (x,y)$, there is a morphism generator $(f,s)\colon (X,x) \to (Y,y)$.
  \end{itemize}
  We want to prove that $\hat{M}$ and $\GrothendieckS{M}$ are isomorphic.
  
  On objects, the proof of bijectivity is as in the previous case.
  We can then define a symmetric monoidal functor 
  $\hat{M} \to \GrothendieckS{M}$ that is 
  bijective on objects, and sends a generating morphism $(f,s)\colon (X,x) \to (Y,y)$ of
  $\hat{M}$ to the morphism $(f,s)\colon (X,x) \to (Y,y)$ in $\GrothendieckS{M}$.
  
  We want to prove that this functor is full and faithful. 
  Faithfulness is again straightforward. 
  Suppose given $(f_1,s_1), (f_2,s_2)\colon (X,x) \to (Y,y)$ in $\hat{M}$. 
  By construction, $(f_1,s_1) = (f_2, s_2)$ in $\GrothendieckS{M}$
  if and only if $f_1 = f_2$ in $\Free{M}$ and $s_1 = s_2$. But this 
  means that $(f_1,s_1) = (f_2,s_2)$ also in $\hat{M}$.
  
  For fullness, take a morphism $(f,s)\colon (X,x) \to (Z,z)$ in 
  $\GrothendieckS{M}$, and notice the following:
  \begin{itemize}
    \item Each morphism $(f,s)$ such that $f$ is a generator, an identity or a symmetry 
    in $\Free{M}$ is also in $\hat{M}$;
    \item If $(f,s)\colon (X,x) \to (Z,z)$ is such that $f = g\cp h$, where $g\colon X \to Y$ 
    and $h\colon Y \to Z$, then by definition of composition in $\Span$, 
    we have $s = (t,u)$ for some $t,u$ with $\Fun{M}g(t)=(x,y)$ and 
    $\Fun{M}h(u)=(y,z)$. This means that $(g,t)\colon (X,x) \to (Y,y)$ and 
    $(h,u)\colon (Y,y) \to (Z,z)$ are morphisms in $\GrothendieckS{M}$;
    \item If $(f,s)\colon (X_1 \Tensor X_2,(x_1, x_2)) \to (Y_1 \Tensor Y_2, (y_1, y_2))$ 
    is such that $f=f_1 \Tensor f_2$, where $f_1\colon X_1 \to Y_1$ and $f_2\colon X_2 \to Y_2$, 
    then $s = (s_1,s_2)$ for some $s_1, s_2$, and we have :
    \begin{equation*}
      ((x_1,x_2),(y_1,y_2))= \Fun{M}f (s)
        = \Fun{M}(f_1 \Tensor f_2) (s_1, s_2)
        = (\Fun{M}f_1 \times \Fun{M}f_2)(s_1, s_2)
    \end{equation*}
    This means that $\Fun{M}f_1(s_1) = (x_1,y_1)$ and 
    $\Fun{M}f_2(s_2) = (x_2,y_2)$, and hence that
    $(f_1,s_1)\colon (X_1,x_1) \to (Y_1,y_1)$ and 
    $(f_2, s_2)\colon (X_2,x_2) \to (Y_2,y_2)$ are morphisms in $\GrothendieckS{M}$.
  \end{itemize}
  By definition, since $\Free{M}$ is free, any morphism $f$ can 
  be decomposed into a composition of monoidal products of morphism generators, 
  symmetries and identities. The points above prove that each 
  of such morphisms is also in $\GrothendieckS{M}$, and hence in 
  $\hat{M}$. So $f$ in $\GrothendieckS{M}$ is the image of $f$ 
  in $\hat{M}$, and the functor is full.
  
  Since our correspondence is bijective on objects 
  and fully faithful, this proves that $\hat{M}$ and $\GrothendieckS{M}$
  are isomorphic as categories. Since $\hat{M}$ is free so is $\GrothendieckS{M}$, 
  completing the proof.
\end{proof}
\begingroup
\def\thelemma{\ref{lem: strict monoidal functor pi SetS}}
\begin{lemma}
  In the case of both $\SetS$ and $\Span$, there is a strict 
  monoidal functor $\pi_M\colon \GrothendieckS{M} \to \Free{M}$ 
  sending a object $(X,x)$ to $X$ and a morphism 
  $f\colon(X,x) \to (Y,y)$ to $f\colon X \to Y$
  (resp.\ $(f,s)\colon(X,x) \to (Y,y)$ to $f\colon X \to Y$).
\end{lemma}
\addtocounter{lemma}{-1}
\endgroup
\begin{proof}
  Let us prove the statement in the case of $\SetS$, $\Span$ being
  analogous.
  In $\GrothendieckS{M}$, the identity on $(X,x)$ is $\Id{X}: (X,x) \to (X,x)$, which
  obviously gets sent to $\Id{X}:X \to X$ in $\Free{M}$.
  Similarly, the monoidal unit $(\TensorUnit, *)$ of $\GrothendieckS{M}$
  is sent to $\TensorUnit$, which is the monoidal unit of $\Free{M}$.
  
  Compositon and tensor product
  are defined by freeness both in $\GrothendieckS{M}$ and $\Free{M}$, 
  and since if $f$ is a generator in $\GrothendieckS{M}$ then it 
  is a generator in $\Free{M}$, we have that $f;g$ (resp. $f \Tensor g$)
  in $\GrothendieckS{M}$ is sent to $f;g$ (resp. $f \Tensor g$) in
  $\Free{M}$. Strictness follows by definition, concluding the proof.
\end{proof}
\begingroup
\def\thetheorem{\ref{thm: reachability is preserved}}
\begin{theorem}
  Let $\NetSem{N}$ be a guarded Petri net (with side effects).
  $(Y,y)$ is reachable from $(X,x)$ if and only if $(Y,y)$ is reachable
  from $(X,x)$ in the net $\UnFree{\GrothendieckS{N}}$.
\end{theorem}
\addtocounter{theorem}{-1}
\endgroup
\begin{proof}
  By definition $(Y,y)$ is reachable from $(X,x)$
  if and only if there is a morphism $f:X \to Y$ in $\Free{N}$ 
  such that $\Fun{N}f(x) = y$ (resp. $\Fun{N}f(s) = (x,y)$ for some $s \in S_f$). 
  Again by definition, this means that
  $f:(X,x) \to (Y,y)$ (resp. $f_s:(X,x) \to (Y,y)$) is a morphism of $\GrothendieckS{N}$.
  Since $\GrothendieckS{N}$ is free, $f$ (resp. $f_s$) can be decomposed
  as a composition of monoidal products of generating morphisms.
  But every generating morphism of $\GrothendieckS{N}$
  corresponds to a transition of $\UnFree{\GrothendieckS{N}}$,
  from which the thesis follows.
\end{proof}
\begingroup
\def\thelemma{\ref{lem: lift F}}
\begin{lemma}
  Let $F\colon \NetSem{M} \to \NetSem{N}$ be a morphism in $\PetriSetS$ (resp.\ in $\PetriSpan$). Then 
  it lifts to strict monoidal functor
  $\LiftSetS{F}\colon \GrothendieckS{M} \to \GrothendieckS{N}$ 
  (resp.\
  $\LiftSpan{F}\colon \GrothendieckS{M} \to \GrothendieckS{N}$), 
  such that the following diagram on the left (resp. on the right) commutes:
  \begin{equation*}
    \begin{tikzpicture}[node distance=1.3cm,>=stealth',bend angle=45,auto]
      \node (1) at (0,1.5) {$\GrothendieckS{M}$};
      \node (2) at (0,0) {$\Free{M}$};
      \node (3) at (3,1.5) {$\GrothendieckS{N}$};
      \node (4) at (3,0) {$\Free{N}$};
      \node (5) at (1.5, -1.5) {$\SetS$};
      \draw[->] (1)--(2) node [midway,left] {$\pi_M$};
      \draw[->] (3)--(4) node [midway,right] {$\pi_N$};
      \draw[->] (2)--(4) node [midway,above] {$F$};
      \draw[->, dashed] (1)--(3) node [midway,above] {$\LiftSetS{F}$};
      \draw[->] (2)--(5) node [midway,left] {$\Fun{M}$};
      \draw[->] (4)--(5) node [midway,right] {$\Fun{N}$};
    \end{tikzpicture}
    \hspace{6em}
    \begin{tikzpicture}[node distance=1.3cm,>=stealth',bend angle=45,auto]
      \node (1) at (0,1.5) {$\GrothendieckS{M}$};
      \node (2) at (0,0) {$\Free{M}$};
      \node (3) at (3,1.5) {$\GrothendieckS{N}$};
      \node (4) at (3,0) {$\Free{N}$};
      \node (5) at (1.5, -1.5) {$\Span$};
      \draw[->] (1)--(2) node [midway,left] {$\pi_M$};
      \draw[->] (3)--(4) node [midway,right] {$\pi_N$};
      \draw[->] (2)--(4) node [midway,above] {$F$};
      \draw[->, dashed] (1)--(3) node [midway,above] {$\LiftSpan{F}$};
      \draw[->] (2)--(5) node [midway,left] {$\Fun{M}$};
      \draw[->] (4)--(5) node [midway,right] {$\Fun{N}$};
    \end{tikzpicture}
  \end{equation*}
\end{lemma}
\addtocounter{lemma}{-1}
\endgroup
\begin{proof}
  We begin with the $\SetS$ case. If $(X,x)$ is in $\GrothendieckS{M}$, then $x \in \Fun{M}X = 
 \Fun{N}(FX)$, so we can define $\LiftSetS{F}(X,x)\coloneqq(FX, x)$ on objects. Clearly it has $\LiftSetS{F}(\TensorUnit_M, *) = 
  (F \TensorUnit_M, *) = (\TensorUnit_N,*)$, and
  \begin{align*}
    \LiftSetS{F}\left((X_1,x_1)\otimes (X_2,x_2)\right)
    &= \LiftSetS{F}\left(X_1\otimes X_2, (x_1, x_2)\right)\\
    &= \left(FX_1\otimes FX_2,(x_1, x_2) \right)\\
    &= (FX_1, x_1)\otimes (FX_2,x_2).\\
  \end{align*}
  So $\LiftSetS{F}$ acts monoidally on objects. On morphisms, 
  we send $f\colon (X,x) \to (Y,y)$ to $\LiftSetS{F}(f)\coloneqq Ff\colon (FX, x) \to (FY, y)$.
  Clearly we have $y = \Fun{M}f(x) = 
  \Fun{N}(Ff)(x)$, proving the mapping is well defined. Functoriality 
  and commutativity of the diagram follow from the definitions.

  \medskip
  \noindent
  For $\LiftSpan{F}$, on objects the proof is as in 
  the previous case. On morphisms, 
  we send $(f,s)\colon (X,x) \to (Y,y)$ to $(Ff,s)\colon (FX, x) \to (FY, y)$.
  Clearly we have $(x,y) = \Fun{M}f(s) = 
  \Fun{N}(Ff)(s)$, proving the mapping is well defined. Functoriality 
  and commutativity of the diagram follow from the definitions.
\end{proof}
\begingroup
\def\thelemma{\ref{lem: hatF SetS is full and faithful}}
\begin{lemma}
  For any map $F\colon(M,\Fun{M})\to(N,\Fun{N})$ in $\PetriSetS$ 
  (respectively in $\PetriSpan$), the functor $F$ is faithful if and only 
  if $\LiftSetS{F}$ is faithful (resp.\ $\LiftSpan{F}$ is faithful). 
  If $F$ is full, then so is $\LiftSetS{F}$ (resp.\ $\LiftSpan{F}$).
\end{lemma}
\addtocounter{lemma}{-1}
\endgroup
\begin{proof}
  We start from $\LiftSetS{F}$.
  For faithfulness, consider 
  $f, g\colon (X,x) \to (Y,y)$. By definition, $\LiftSetS{F}f=Ff$ and $\LiftSetS{F}g = Fg$,
  and $\LiftSetS{F}$ is faithful if and only if $F$ is. 
  
  Now suppose $F$ is full, and consider $f\colon (FX, x) \to (FY,y)$ in 
  $\GrothendieckS{N}$. Then the corresponding $f$ in
  $\Free{N}$ is a morphism $FX \to FY$, and since $F$ 
  is full, there is a $g\colon X \to Y$ in $\Free{M}$ such that 
  $Fg = f$. Moreover, we have that $y = \Fun{N}f(x)
  = \Fun{N}(Fg)(x)  = \Fun{M}g(x)$, 
  so there is a morphism $g\colon (X,x) \to (Y,y)$ in $\GrothendieckS{M}$.
  Clearly it is $\LiftSetS{F}g = Fg = f$, proving that $\LiftSetS{F}$ is full.

  \bigskip
  \noindent
  Now we focus on $\LiftSpan{F}$. For faithfulness, suppose given
  $(f,s),(g,t)\colon (X,x) \to (Y,y)$; they are equal iff $f=g$ and $s=t$. By definition, $\hat{F}(f,s)=(Ff,s)$ and
  $\hat{F}(g,t) = (Fg,t)$,
  and hence $\hat{F}(f,s) = \hat{F}(g,t)$ iff $Ff=Fg$ in $\Free{N}$ and $s=t$. Thus $F$ is faithful iff $\hat{F}$ is faithful.
  
    Now suppose $F$ is full, and consider $(f,s)\colon (FX, x) \to (FY,y)$ in 
  $\GrothendieckS{N}$. Then the corresponding $f$ in
  $\Free{N}$ is a morphism $FX \to FY$, and since $F$ 
  is full, there is a $g: X \to Y$ in $\Free{M}$ such that 
  $Fg = f$. Moreover, we have that $(x,y) = \Fun{N}f(s)
  = \Fun{N}(Fg)(s) = \Fun{M}g(s)$, 
  so there is a morphism $(g,s)\colon (X,x) \to (Y,y)$ in $\GrothendieckS{M}$.
  Clearly it is $\hat{F}(g,s) = (Fg,s) = (f,s)$, proving that $\hat{F}$ is full and concluding 
  the proof.
\end{proof}
\begingroup
\def\thelemma{\ref{lem: hat F SetS is transition-preserving}}
\begin{lemma}
  If $F$ is transition-preserving, so are $\LiftSetS{F}$ and $\LiftSpan{F}$.
\end{lemma}
\addtocounter{lemma}{-1}
\endgroup
\begin{proof}
  The proof is obvious considering that, by construction, $f\colon (X,x) \to (Y,y)$ 
  is a generator (resp. a symmetry) in $\GrothendieckS{N}$ if 
  and only if it is a generator (resp. a symmetry) in $\Free{N}$.

  An analogous argument holds for $\LiftSpan{F}$.
\end{proof}
\begingroup
\def\thelemma{\ref{lem: hat F SetS is injective on objects}}
\begin{lemma}
  If $F$ is injective on objects, so are $\LiftSetS{F}$ and $\LiftSpan{F}$.
\end{lemma}
\addtocounter{lemma}{-1}
\endgroup
\begin{proof}
  Suppose $\LiftSetS{F}(X,x) = \hat{F}(Y,y)$. By definition, 
  this means $(FX,x) = (FY,y)$, from which it follows that $x = y$ and $FX = FY$. Since $F$ is 
  injective on objects, then $X = Y$ and so $(X,x) = (Y,y)$, proving that $\LiftSetS{F}$ is injective 
  on objects as well.

  The span case is analogous.
\end{proof}
\begingroup
\def\thelemma{\ref{lem: hat F preserves identifications}}
\begin{lemma}
  If $\NetSem{N}$ is an identification of $\NetSem{M}$ via
  $F$ and witnesses $O, l, r$, then $\GrothendieckS{N}$ is an 
  identification of $\GrothendieckS{M}$ via $\LiftSetS{F}$ and 
  witnesses $\UnFree{\Grothendieck{(l \cp \Fun{M})}}$, $\LiftSetS{l}, \LiftSetS{r}$.
  The span case is analogous.
\end{lemma}
\addtocounter{lemma}{-1}
\endgroup
\begin{proof}
  If $\NetSem{N}$ is an identification of $\NetSem{M}$ via
  $F$ and witnesses $O, l, r$, then by definition:
  \begin{itemize}
    \item $l$ and $r$ are transition-preserving;
    \item $l \cp \Fun{M} = r \cp \Fun{M}$;
    \item $F:\Free{M} \to \Free{N}$ is the coequalizer of $l,r: \Free{O} \to \Free{M}$.
  \end{itemize}
  From the first point we get that $\GrothendieckS{l \cp M} = \GrothendieckS{r \cp M}$,
  so we can consider:
  \begin{equation*}
    \scalebox{0.8}{
      \begin{tikzpicture}[node distance=2cm,>=stealth',bend angle=45,auto]
        \node (1) at (0,0) {$\Free{M}$};
        \node (2) [right of=1] {$\Free{N}$};
        \node (w) [left of=1] {$\Free{O}$};

        \node[red!60!black](1a) [above of=1] {$\GrothendieckS{M}$};
        \node[red!60!black](2a) [above of=2] {$\GrothendieckS{N}$};
        \node[red!60!black](wa) [above of=w] {$\GrothendieckS{l \cp M}$};

        \node(1b) [below of=1] {$\Semantics$};

        \draw[transform canvas={yshift=0.5ex},->] (w) to node {$l$} (1);
        \draw[transform canvas={yshift=-0.5ex},->](w) to node[swap] {$r$} (1); 
        \draw[red!60!black, transform canvas={yshift=0.5ex},->] (wa) to node {$\LiftSetS{l}$} (1a);
        \draw[red!60!black, transform canvas={yshift=-0.5ex},->](wa) to node[swap] {$\LiftSetS{r}$} (1a); 

        \draw[dotted, ->] (1) to node {$F$} (2);
        \draw[red!60!black, dotted, ->] (1a) to node {$\LiftSetS{F}$} (2a);

        \draw[red!60!black, ->] (1a) to node {$\pi$} (1);
        \draw[red!60!black, ->] (2a) to node {$\pi$} (2);
        \draw[red!60!black, ->] (wa) to node {$\pi$} (w);

        \draw[->] (w) to node[swap] {$l \cp \Fun{M}$} (1b);
        \draw[->] (1) to node[swap] {$\Fun{M}$} (1b);
        \draw[dashed, ->] (2) to node {$\Fun{N}$} (1b);
      \end{tikzpicture}
   }
  \end{equation*}
  We need to prove that $\LiftSetS{F}:\Free{M} \to \Free{N}$ is the coequalizer of 
  $\LiftSetS{l}$ and $\LiftSetS{r}$. 
  Note that the following proof would have been much simpler
  if we had an adjunction between Petri nets and FSSMCs.
  Unfortunately we do not have such an adjunction as pointed
  out in~\cite{Genovese2019c}.

  According to~\cite{Genovese2019}, Lemma 4.2, generating 
  objects of $\NetSem{N}$ can be written as
  $\GObj{\Free{M}}/\simeq_o$, where $\simeq_o$
  is the equivalence relation generated by 
  \begin{equation*}
    \forall X \in \GObj{\Free{O}}.(lX = rX)
  \end{equation*}
  Now,  objects of $\GrothendieckS{M}$ are pairs $(X,x)$ with
  $X$ object of $\Free{M}$ and $x \in \Fun{M}X$. Consider then 
  $\GObj{\GrothendieckS{M}/\simeq_{o'}}$, where $\simeq_{o'}$
  is the equivalence relation generated by 
  \begin{equation*}
    \forall (X,x) \in \GObj{(\GrothendieckS{l \cp M})}.(\LiftSetS{l}(X,x) = \LiftSetS{r}(X,x))
  \end{equation*}
  Unpacking definitions, $\LiftSetS{l}(X,x) = \LiftSetS{r}(X,y)$ iff $lX = rX$ and $x = y$,
  So for each generating object $(X,x)$ of $\GrothendieckS{M}$ we have 
  $[(X,x)]_{\simeq_{o'}} = ([X]_{\simeq_o}, x)$. This proves that
  \begin{equation*}
    \GObj{\GrothendieckS{M}/\simeq_{o'}} = \GObj{\GrothendieckS{N}}
  \end{equation*}
  Again according to~\cite{Genovese2019}, Lemma 4.2,
  given a generating morphism $f$ of $\Free{O}$, and 
  denoting with $f_l$ and $f_r$ the generating morphisms 
  of $\Free{M}$ such that $lf = \sigma \cp f_l \cp \sigma'$ and 
  $rf = \varsigma \cp f_r \cp \varsigma'$, generating 
  morphisms of $\Free{N}$ are 
  $\GMor{\Free{M}}/\simeq_m$, where $\simeq_m$
  is the equivalence relation generated by 
  \begin{equation*}
  \forall f \in \GMor{\Free{O}}.(f_l = f_r)
  \end{equation*}
  Since $\LiftSetS{l}f$ (resp. $\LiftSetS{r}f$) agrees with $lf$ whenever 
  $f$ exists in $\GrothendieckS{l \cp M}$ (resp. in $\GrothendieckS{r \cp M}$),
  requiring  $\simeq_{m'}$ to be the equivalence relation 
  on $\GMor{\GrothendieckS{M}}$ defined by 
  \begin{equation*}
    \forall f \in \GMor{(\GrothendieckS{l \cp M})}.(f_l = f_r)
  \end{equation*}
  We have that, whenever they exist, $f \simeq_{m'} g$ in $\GMor{\GrothendieckS{M}}$ iff
  $f \simeq_m g $ in $\GMor{\Free{M}}$. Considering that $f$ exists in $\GMor{\GrothendieckS{M}}$
  if and only if $[f]_m$ exists in $\GMor{\GrothendieckS{N}}$, we have proven that:
  \begin{equation*}
    \GMor{\GrothendieckS{M}}/\simeq_{m'} = \GMor{\GrothendieckS{N}}
  \end{equation*}
  Applying one last time the characterization of~\cite{Genovese2019}, Lemma 4.2,
  We have that $\LiftSetS{F}:\GrothendieckS{M} \to \GrothendieckS{N}$ is the coequalizer
  of $\LiftSetS{l}$ and $\LiftSetS{r}$.
  
  In $\Span$ the proof on objects in analogous. On morphisms, notice that 
  since $\LiftSetS{l}(f,s) = (lf,s)$ (resp. $\LiftSetS{r}(f,s) = (rf,s)$ )
  requiring  $\simeq_{m'}$ to be the equivalence relation 
  on $\GMor{\GrothendieckS{M}}$ defined by 
  \begin{equation*}
    \forall f \in \GMor{(\GrothendieckS{l \cp M})}.((f_l, s) = (f_r, s))
  \end{equation*}
  We have that, whenever they exist, $(f,s) \simeq_{m'} (g,s)$ in $\GMor{\GrothendieckS{M}}$ iff
  $f \simeq_m g $ in $\GMor{\Free{M}}$. Considering that $(f,s)$ exists in $\GMor{\GrothendieckS{M}}$
  if and only if $[f]_m$ exists in $\GMor{\GrothendieckS{N}}$, we have proven that:
  \begin{equation*}
    \GMor{\GrothendieckS{M}}/\simeq_{m'} = \GMor{\GrothendieckS{N}}
  \end{equation*}
  The rest of the proof proceeds as in the $\SetS$ case.
\end{proof}
\begingroup
\def\thelemma{\ref{lem: hat F preserves addition of generating morphisms}}
\begin{lemma}
  Let $\NetSem{M}$  be an addition of generating
  morphisms to $\NetSem{K}$ via synchronization witness
  $W, w$. 
  Then $\GrothendieckS{M}$ is an addition of 
  generating morphisms to $\GrothendieckS{K}$ via 
  synchronization witness $\UnFree{\Grothendieck{(w \cp \Fun{K})}},
  \LiftSetS{w}$. The span case is analogous.
\end{lemma}
\addtocounter{lemma}{-1}
\endgroup
\begin{proof}
 To prove this, consider the following commutative diagram,
 where we use the notation developed in~\cite{Genovese2019}:
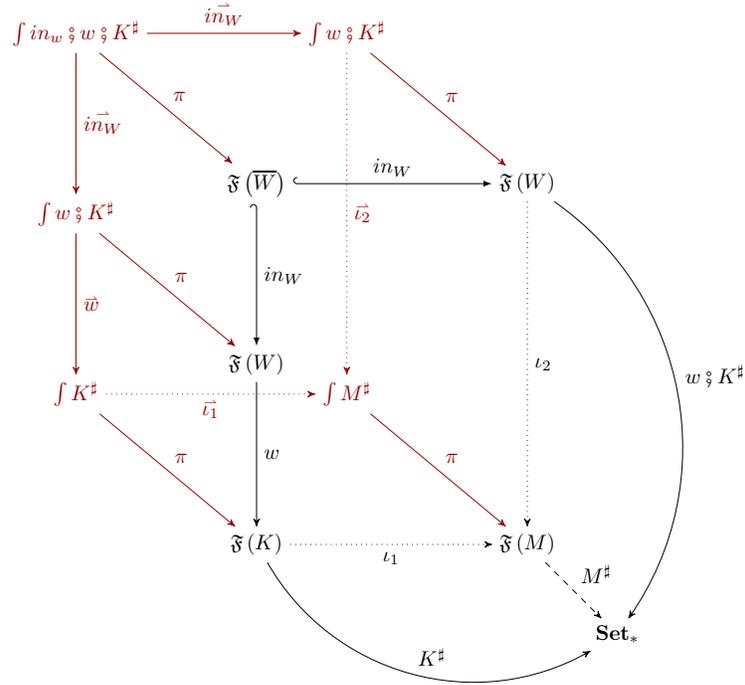
 \begin{equation*}
  \scalebox{0.8}{
    \begin{tikzpicture}[node distance=2cm,>=stealth',bend angle=45,auto]
      \node (in1) at (0,6) {$\Free{\overline{W}}$};
      \node (in2) at (4.5,6) {$\Free{W}$};
      \node (in3) at (0,3) {$\Free{W}$};
      \node (in4) at (0,0) {$\Free{K}$};
      \node (in5) at (4.5,0) {$\Free{M}$};
  
      \node (sem) at (6, -1.5) {$\SetS$};
  
      \node[red!60!black] (out1) at (-3,8.5) {$\Grothendieck{in_w \cp w \cp \Fun{K}}$};
      \node[red!60!black] (out2) at (1.5,8.5) {$\Grothendieck{w \cp \Fun{K}}$};
      \node[red!60!black] (out3) at (-3,5.5) {$\Grothendieck{w \cp \Fun{K}}$};
      \node[red!60!black] (out4) at (-3,2.5) {$\Grothendieck{\Fun{K}}$};
      \node[red!60!black] (out5) at (1.5,2.5) {$\Grothendieck{\Fun{M}}$};
  
      \draw[right hook-latex] (in1) to node {$in_W$} (in2);
      \draw[left hook-latex] (in1) to node {$in_W$} (in3);
      \draw[->] (in3) to node {$w$} (in4);
      \draw[->, dotted] (in4) to node[below] {$\iota_1$} (in5);
      \draw[->, dotted] (in2) to node {$\iota_2$} (in5);
      \draw[->, dashed] (in5) to node {$\Fun{M}$} (sem);
      \draw[->, bend right] (in4) to node {$\Fun{K}$} (sem);
      \draw[->, bend left] (in2) to node {$w \cp \Fun{K}$} (sem);
  
      \draw[->, red!60!black] (out1) to node {$\pi$} (in1);
      \draw[->, red!60!black] (out2) to node {$\pi$} (in2);
      \draw[->, red!60!black] (out3) to node {$\pi$} (in3);
      \draw[->, red!60!black] (out4) to node {$\pi$} (in4);
      \draw[->, red!60!black] (out5) to node {$\pi$} (in5);
  
      \draw[->, red!60!black] (out1) to node {$\LiftSetS{in_W}$} (out2);
      \draw[->, red!60!black] (out1) to node {$\LiftSetS{in_W}$} (out3);
      \draw[->, red!60!black] (out3) to node {$\LiftSetS{w}$} (out4);
      \draw[->, dotted, red!60!black] (out4) to node[below] {$\LiftSetS{\iota_1}$} (out5);
      \draw[->, dotted, red!60!black] (out2) to node {$\LiftSetS{\iota_2}$} (out5);

    \end{tikzpicture}
  }
 \end{equation*}
 The dotted arrows are a pushout, and can be characterized 
 as identifications on the coproduct. 
 We proved in~\ref{lem: hat F preserves identifications} that identifications
 are preserved, while we defer the proof that coproducts are preserved 
 in the proof of Threorem~\ref{thm: mapping to Term}.
 Hence, $\GrothendieckS{M}$ is the pushout of
 $\LiftSetS{in_W}$ and $\LiftSetS{in_W} \cp \LiftSetS{w}$.
 The claim follows immediately by noticing that:
 \begin{equation*}
  \Grothendieck{in_W \cp w \cp \Fun{K}} = \overline{\Grothendieck{w \cp \Fun{K}}}
  \qquad
  \LiftSetS{in_W} = in_{{w \cp \Fun{K}}}
 \end{equation*}
 which is obvious since $\Free{\overline{W}}$ and $\Free{W}$
 coincide on generating objects and $\Free{\overline{W}}$ has
 no generating morphisms.
\end{proof}
\begingroup
\def\thelemma{\ref{lem: hat F preserves erasing of generating morphisms}}
\begin{lemma}
  Let $\NetSem{K}$ be an erasing of generating
  morphisms from $\NetSem{N}$ via a subnet $N_w$.
  Then $\GrothendieckS{K}$ is an erasing 
  of generators from $\GrothendieckS{N}$ via
  $\Grothendieck sub_{N_w} \cp \Fun{N}$. The span case is 
  analogous.
\end{lemma}
\addtocounter{lemma}{-1}
\endgroup
\begin{proof}
  The proof is easy by considering the following diagram:
  \begin{equation*}
    \scalebox{0.8}{
      \begin{tikzpicture}[node distance=2cm,>=stealth',bend angle=45,auto]
        \node (in1) at (0,6) {$\Free{N_w}$};
        \node (in2) at (4.5,6) {$\Free{\overline{N_w}}$};
        \node (in3) at (0,0) {$\Free{N}$};
        \node (in4) at (4.5,0) {$\Free{K}$};

        \node (sem) at (-1.5, -1.5) {$\SetS$};

        \node[red!60!black] (out1) at (3,8.5) {$\Grothendieck{sub_{N_w} \cp \Fun{N}}$};
        \node[red!60!black] (out2) at (7.5,8.5) {$\Grothendieck{in_{N_w} \cp sub_{N_w} \cp \Fun{N}}$};
        \node[red!60!black] (out3) at (3,2.5) {$\Grothendieck{\Fun{N}}$};
        \node[red!60!black] (out4) at (7.5,2.5) {$\Grothendieck{\Fun{K}}$};

        \draw[left hook-latex] (in2) to node[swap] {$in_{N_w}$} (in1);
        \draw[left hook-latex, dotted] (in1) to node {$sub_{N_w}$} (in3);
        \draw[left hook-latex, dotted] (in4) to node[below] {$sub_K$} (in3);
        \draw[->] (in2) to node {$k$} (in4);
        \draw[->, bend left] (in4) to node {$\Fun{K}$} (sem);
        \draw[->, dashed] (in3) to node {$\Fun{N}$} (sem);
        \draw[->, bend right] (in1) to node[swap] {$sub_{N_w} \cp \Fun{N}$} (sem);

        \draw[->, red!60!black] (out1) to node {$\pi$} (in1);
        \draw[->, red!60!black] (out2) to node {$\pi$} (in2);
        \draw[->, red!60!black] (out3) to node {$\pi$} (in3);
        \draw[->, red!60!black] (out4) to node {$\pi$} (in5);

        \draw[->, red!60!black] (out2) to node[swap] {$\LiftSetS{in_{N_w}}$} (out1);
        \draw[left hook-latex, dotted, red!60!black] (out1) to node {$\LiftSetS{sub_{N_w}}$} (out3);
        \draw[left hook-latex, dotted, red!60!black] (out4) to node[below] {$\LiftSetS{sub_K}$} (out3);
        \draw[->, red!60!black] (out2) to node {$\LiftSetS{k}$} (out4);

      \end{tikzpicture}
    }
  \end{equation*}
 Again, the black arrows are the definition
 of erasing of generators, and are a pushout square.
 This can be characterized as the coequalizer or a coproduct,
 where all the morphisms involved are transition-preserving.
 The maroon objects and arrows are given by the
 Grothendieck construction and preserve the pushout.
 The thesis follows by noticing that:
 \begin{align*}
  \Grothendieck{in_{N_w} \cp sub_{N_w} \cp \Fun{N}} &= 
  \overline{\Grothendieck{sub_{N_w} \cp \Fun{N}}}\\
  \LiftSetS{in_{N_w}} &= in_{{sub_{N_w} \cp \Fun{N}}}\\
  \LiftSetS{sub_{N_w}} &= sub_{\Grothendieck sub_{N_w} \cp \Fun{N}}
  \qedhere
 \end{align*}
\end{proof}
\begingroup
\def\thetheorem{\ref{thm: mapping to Term}}
\begin{theorem}
  Denote with $\TermCat$ the terminal category, together with the
   trivial symmetric monoidal structure on it. There is a faithful, 
   strong monoidal functor $\EmbSpan:\PetriSetS \to \PetriTerm$
   defined as follows:
  \begin{itemize}
    \item On objects, it sends $\NetSem{M}$ to 
    $\NetSemTerm{\UnFree{\GrothendieckS{M}}}$.
    \item On morphisms, we send the functor
      $F:\NetSem{M} \to \NetSem{N}$ to the 
      functor\footnote{
        To be absolutely precise, we are referring to
        the functor $\Free{\UnFree{\GrothendieckS{M}}} \simeq 
        \GrothendieckS{M} \xrightarrow{\LiftSetS{F}} \GrothendieckS{N} 
        \simeq \Free{\UnFree{\GrothendieckS{N}}}$.}
      \begin{equation*}
        \LiftSetS{F}: \NetSemTerm{\UnFree{\GrothendieckS{M}}} 
          \to  \NetSemTerm{\UnFree{\GrothendieckS{N}}}
      \end{equation*}
  \end{itemize}
  Similarly, there is a faithful, 
  strong monoidal functor
  $\EmbSpan:\PetriSpan \to \PetriTerm$
  defined as follows:
  \begin{itemize}
    \item On objects, it sends $\NetSem{M}$ to 
    $\NetSemTerm{\UnFree{\GrothendieckS{M}}}$.
    \item On morphisms, we send the functor
      $F:\NetSem{M} \to \NetSem{N}$ to the 
      functor\footnote{
        To be absolutely precise, we are referring to
        the functor $\Free{\UnFree{\GrothendieckS{M}}} \simeq 
        \GrothendieckS{M} \xrightarrow{\LiftSpan{F}} \GrothendieckS{N} 
        \simeq \Free{\UnFree{\GrothendieckS{N}}}$.}
      \begin{equation*}
        \LiftSpan{F}: \NetSemTerm{\UnFree{\GrothendieckS{M}}} 
          \to  \NetSemTerm{\UnFree{\GrothendieckS{N}}}.
      \end{equation*}
  \end{itemize}
\end{theorem}
\addtocounter{theorem}{-1}
\endgroup
\begin{proof}
  The proofs for $\SetS$ and $\Span$ are very similar,
  so we just provide the one for $\SetS$.
  Clearly if $F$ is the identity functor $\NetSem{M} \to \NetSem{M}$
  then so is $\LiftSetS{F}$. For composition, 
  consider $F\colon\NetSem{M} \to \NetSem{M}{\Fun{N}}$ and 
  $G\colon \NetSem{N} \to \NetSem{P}$. We have to prove that $\LiftSetS{F\cp G} = 
  \LiftSetS{F}\cp\LiftSetS{G}$. On objects, $\LiftSetS{F\cp G}$ sends 
  $(X,x)$ in $\GrothendieckS{M}$ to $((F\cp G)X, x)$ in $\GrothendieckS{P}$, 
  so it coincides with $\LiftSetS{F}\cp\LiftSetS{G}$. Now consider a morphism $f\colon (X,x) \to (Y,y)$
  in $\GrothendieckS{M}$. This is sent by $\LiftSetS{F}$ to $Ff\colon (FX,x) \to (FY,y)$, and 
  applying $\LiftSetS{G}$ to it one gets $G(Ff)\colon (G(FX),x) \to (G(FY), y)$. Since 
  $G(F(\_))$ is $(F\cp G)(\_)$, we are done. This proves that $\EmbSetS$
  is a functor. Faithfulness is trivial.

  Now we focus on monoidality. First of all we have to prove that 
  \begin{equation*}
    \EmbSetS\left(\NetSem{M} \Tensor \NetSem{N}\right) \simeq 
      \EmbSetS\NetSem{M} \Tensor \EmbSetS\NetSem{N}
  \end{equation*}
  Remembering from Definition~\ref{def: PetriS} 
  that for each choice of semantics $\Semantics$ 
  the monoidal structure on $\PetriS{\Semantics}$ is defined 
  in terms of coproduct of symmetric monoidal categories, 
  and hence from the coproduct of the underlying nets, this means that:
  \begin{align*}
    &\NetSemTerm{\UnFree{\Grothendieck{[\Fun{M},\Fun{N}]}}}=\\
    &\hspace{10em}= \EmbSetS(M + N, [\Fun{M},\Fun{N}])\\ 
    &\hspace{10em}= \EmbSetS(\NetSem{M} \Tensor \NetSem{N})\\
    &\hspace{10em}\simeq \EmbSetS\NetSem{M} \Tensor \EmbSetS\NetSem{N}\\
    &\hspace{10em}= \NetSemTerm{\UnFree{\GrothendieckS{M}}} 
      \Tensor \NetSemTerm{\UnFree{\GrothendieckS{N}}}\\
    &\hspace{10em}= \left(\UnFree{\GrothendieckS{M}} + \UnFree{\GrothendieckS{N}}, 
    \left[\Fun{\UnFree{\GrothendieckS{M}}}, \Fun{\UnFree{\GrothendieckS{N}}}\right]\right)
      \end{align*}
  Since $\UnFree{\_}$ preserves isomorphisms and coproducts, it is sufficient to prove:
  \begin{equation*}
    \Grothendieck{[\Fun{M},\Fun{N}]} 
      \simeq \GrothendieckS{M} + \GrothendieckS{N}
  \end{equation*}
  \begin{itemize}
    \item By definition, objects of $\Grothendieck{[\Fun{M},\Fun{N}]}$ 
    are pairs $(X,x)$ with $X \in \Free{M+N} \simeq \Free{M} + \Free{N}$ and 
    $x \in [\Fun{M},\Fun{N}]X$. This is clearly isomorphic to 
    $\Obj{\GrothendieckS{M}} \sqcup \Obj{\GrothendieckS{N}}$.

    \item Again by definition, we have
    \begin{align*}
      &\Hom{\Grothendieck{[\Fun{M},\Fun{N}]}}
        {(X,x)}{(Y,y)} :=\\
      & \qquad := \Suchthat{f \in \Hom{\Free{M}+\Free{N}}
        {X}{Y}}{[\Fun{M},\Fun{N}]f(x) = y}
    \end{align*}    
    This follows noting that by definition the 
    set of morphisms of $\Free{M} + \Free{N}$ is the disjoint 
    union of the sets of morphisms of $\Free{M}$ and $\Free{N}$.
  \end{itemize}
  Then we have to prove that $\EmbSetS(F \Tensor G) = 
  \EmbSetS F \Tensor \EmbSetS G$. Unrolling definitions 
  this amounts to prove that $\LiftSetS{F + G} = \LiftSetS{F} + \LiftSetS{G}$, 
  which is obvious.

  Finally, we need to prove that $\EmbSetS$ preserves the monoidal 
  unit. Notice that for each choice of semantics $\Semantics$ the monoidal 
  unit in $\PetriS{\Semantics}$ is taken to be $(\emptyset, \Fun{\emptyset}{\Semantics})$.
  $\Free{\emptyset}$ is the free category consisting of only 
  the monoidal unit $I$, and $\Fun{\emptyset}{\Semantics}: 
  \Free{\emptyset} \to \Semantics$ sends the monoidal unit 
  to the monoidal unit and its identity to itself. In our case, the monoidal 
  unit of $\PetriSetS$ is $(\emptyset, \Fun{\emptyset})$, 
  with $\Fun{\emptyset}$ sending the monoidal unit 
  of $\Free{\emptyset}$ to the 
  singleton set $\{*\}$ in $\SetS$. In particular, this means that 
  $\GrothendieckS{\emptyset} \simeq \Free{\emptyset}$, proving 
  that
  \begin{equation*}
    \EmbSetS\NetSem{\emptyset} = \NetSemTerm{\UnFree{\GrothendieckS{\emptyset}}}\\
    \simeq \NetSemTerm{\UnFree{\Free{\emptyset}}}\\
    \simeq \NetSemTerm{\emptyset} \qedhere
  \end{equation*} 
\end{proof}
\begingroup
\def\theproposition{\ref{lem: natural transformations from emb to for}}
\begin{proposition}
  Denote with 
  \begin{equation*}
    \ForSetS:\PetriSetS \to \PetriTerm \qquad \ForSpan:
    \PetriSpan \to \PetriTerm
  \end{equation*}
  the functors defined 
  by sending each Petri net $\NetSem{M}$
  to $\NetSemTerm{M}$.
  Then there are natural transformations:
  \begin{equation*}
    \begin{tikzpicture}[node distance=1.3cm,>=stealth',bend angle=45,auto]
      \node (1) at (0,0) {$\PetriSetS$};
      \node (2) at (3,0) {$\PetriTerm$};
      \node[rotate=-90] (3) at (1.5,0) {$\Rightarrow$};
      \node(4) at (1.8,0) {$\pi$};
      \draw[->, bend right] (1) to node [midway,below] {$\ForSetS$} (2);
      \draw[->, bend left] (1) to node [midway,above] {$\EmbSetS$} (2);
    \end{tikzpicture}
    \qquad
    \begin{tikzpicture}[node distance=1.3cm,>=stealth',bend angle=45,auto]
      \node (1) at (0,0) {$\PetriSpan$};
      \node (2) at (3,0) {$\PetriTerm$};
      \node[rotate=-90] (3) at (1.5,0) {$\Rightarrow$};
      \node(4) at (1.8,0) {$\pi$};
      \draw[->, bend right] (1) to node [midway,below] {$\ForSpan$} (2);
      \draw[->, bend left] (1) to node [midway,above] {$\EmbSpan$} (2);
    \end{tikzpicture}
  \end{equation*}
\end{proposition}
\addtocounter{proposition}{-1}
\endgroup
\begin{proof}
  For each object $\NetSem{M}$ in $\PetriSetS$, we set 
  \begin{equation*}
    \pi_{\NetSem{M}}: \NetSemTerm{\UnFree{\GrothendieckS{M}}} 
    \to \NetSemTerm{M}   
  \end{equation*}
  to be the functor $\Free{\UnFree{\GrothendieckS{M}}} 
  \simeq \GrothendieckS{M} \xrightarrow{\pi_M} \Free{M}$, where $\pi_M$ 
  is defined as in Lemma~\ref{lem: lift F}. The naturality condition 
  follows from Lemma~\ref{lem: lift F} as well.
  The span case is analogous.
\end{proof}
\end{document}